\pdfoutput=1
\documentclass[11pt,a4paper,reqno]{article}

\usepackage{amssymb}
\usepackage{amsmath}
\usepackage{amsthm}
\usepackage{booktabs}
\usepackage[dvipsnames]{xcolor}
\usepackage{mathrsfs}
\usepackage[margin=1in]{geometry}
\usepackage{authblk}
\usepackage{enumitem}
\setitemize{itemsep=0em}
\setenumerate{itemsep=0em}
\setenumerate[1]{label=(\arabic*)}
\usepackage{hyperref}
\hypersetup{
  colorlinks = true,
  urlcolor = blue,
  linkcolor = purple,
  citecolor = ForestGreen
}
\usepackage{orcidlink}

\numberwithin{equation}{section}

\theoremstyle{plain}
\newtheorem{theorem}{Theorem}[section]
\newtheorem{corollary}[theorem]{Corollary}
\newtheorem{proposition}[theorem]{Proposition}
\newtheorem{lemma}[theorem]{Lemma}
\theoremstyle{definition}
\newtheorem{definition}[theorem]{Definition}
\newtheorem{example}[theorem]{Example}
\newtheorem{remark}[theorem]{Remark}

\newcommand{\Gal}{\mathrm{Gal}}

\allowdisplaybreaks

\title{Galilean Cartan connections for second order differential equations with a metric}
\author[1,2]{Brad Lackey\orcidlink{0000-0002-3823-8757}}
\affil[1]{Microsoft Quantum, Redmond, WA 98052, USA}
\affil[2]{Department of Physics, University of Washington, Seattle, WA 98195, USA}

\date{\today}

\begin{document}

\maketitle

\begin{abstract}
We study Cartan geometries modeled on the first jet space of curves, a homogeneous space of the Galilean group.
Under a partial integrability condition, such a geometry has laboratory coordinates, in which its geodesics are the solutions of second order differential equations and its connection determines a metric on the spatial directions, the fundamental tensor.
We classify, by their torsion, the connections with prescribed laboratory coordinates, differential equations, and fundamental tensor.
The covariant derivative of one of them, the normal Galilean connection, is the linear connection of Crampin, Mart\'inez, and Sarlet corrected to preserve the fundamental tensor, and one of its torsion components is the Jacobi endomorphism of the equations.
For a second order differential equation on a Riemannian manifold, the normal Galilean connection exists globally, and the Helmholtz conditions for a Lagrangian whose multiplier matrix is the metric become two conditions on its torsion.
A global Lagrangian then exists if and only if a closed two-form on the manifold is exact.
Finally, we derive the Jacobi equation and give conditions under which its solutions are bounded.


\end{abstract}

\section{Introduction}
\label{sec:introduction}

In \cite{lackey2000metric}, motivated by models in biology, we applied Cartan's method of equivalence to a system of second order ordinary differential equations together with a metric, and obtained a preferred connection whose torsion and curvature, with their covariant derivatives, generate the invariants of the two structures.

The geometry of a system of second order ordinary differential equations under the changes of variables \(\bar t=t\), \(\bar x=\bar x(t,x)\), Chern's problem (B) \cite{chern1939geometrie}, is a Cartan geometry modeled on the first jet space, regarded as a homogeneous space of the group \(G_n\) of Section~\ref{sec:orthogonal-cartan}, whose elements change time, position, and velocity affinely and the spatial frame by a general linear transformation \cite{grant2012galilean}.
In \cite{grant2012galilean} we proved that, under a partial integrability condition on the torsion, such a geometry has \emph{laboratory coordinates} \((t,x^i,v^i)\), in which, for a suitable section, the components \(\tau\), \(\omega^i\), \(\phi^i\) of the connection are \(dt\), \(dx^i-v^i\,dt\), and \(dv^i\) modulo \(dt\) and \(dx^j\), and we classified the connections satisfying this condition by their torsion.
That paper treats the differential equations alone, with general linear frames.
Here we use orthogonal frames, so that the geometry also determines a metric, the fundamental tensor; the existence of laboratory coordinates (Theorem~\ref{thm:laboratory}) and the geodesic equations (Proposition~\ref{prop:geodesic-equations}) are proved as in \cite{grant2012galilean}, with \(dx^i-v^i\,dt\) and \(dv^i\) now multiplied by a symmetric positive definite matrix-field.
We classify, together with their torsion, the connections of these geometries that have prescribed laboratory coordinates, geodesic equations, and fundamental tensor; the connection of \cite{lackey2000metric} is one of them (Section~\ref{sec:curvature-bianchi}).
Unlike the equivalence problems of Fels \cite{fels1995equivalence} and of Doubrov, Komrakov, and Morimoto \cite{doubrov1999equivalence}, which classify differential equations up to point and contact transformations, respectively, we fix both the equations and the fundamental tensor.

Our Galilean manifolds are modeled on the \((2n+1)\)-dimensional jet space
of time, positions, and velocities. This differs from, for example, the
\((n+1)\)-dimensional spacetime Galilei structures of
K{\"u}nzle \cite{kunzle1972galilei} and Newton--Cartan geometry.
The latter is related to Lorentzian geometry through Bargmann
structures \cite{duval1985bargmann} and the gauging of the Bargmann
algebra \cite{andringa2011newtonian}.

We call a Galilean connection \emph{regular} when its torsion satisfies the partial integrability condition under which laboratory coordinates exist (Theorem~\ref{thm:laboratory}). Our main results are the following.
The fundamental theorem (Theorem~\ref{thm:metric-general}) classifies the regular Galilean connections having prescribed coordinates as laboratory coordinates, prescribed geodesic equations, and a prescribed fundamental tensor: metric compatibility determines the symmetric parts of the torsion coefficients \(B\), \(M\), and \(P\), and the connections are parameterized by the torsion coefficient \(D\) together with the skew parts of \(B\), \(M\), and \(P\). Exactly one of them, the normal Galilean connection, has \(D=0\) and \(B\), \(M\), \(P\) symmetric in their first two indices (Corollary~\ref{cor:local-normal}).
Its covariant derivative is the linear connection that Crampin, Mart\'inez, and Sarlet attach to the differential equations \cite{crampin1996linear}, corrected by the torsion coefficients \(B\), \(P\), and \(-M\) so that it preserves the fundamental tensor; in particular, \(B\) is half the dynamical covariant derivative of the fundamental tensor (Section~\ref{sec:curvature-bianchi}).
For a Riemannian manifold \((M,g)\), every second order differential equation on \(M\) parameterized by \(t\) is the vector field of a unique normal Galilean connection on \(J^1(\mathbb R,M)\) with fundamental tensor \(g\) (Proposition~\ref{prop:global-normal}).
For Lagrangians with multiplier matrix \(g\), the Helmholtz conditions of the inverse problem \cite{douglas1941solution,sarlet1982helmholtz} become conditions on the torsion, as they do for the connection of \cite[\S4]{lackey2000metric}: such an equation is the Euler--Lagrange equation of a Lagrangian with multiplier matrix \(g\) near each point if and only if \(B=0\) and the Jacobi endomorphism \(\mathcal R\) is self-adjoint (Proposition~\ref{prop:local-inverse}).
The obstruction of Takens and of Anderson and Duchamp to a global Lagrangian \cite{takens1979global,anderson1980existence} then vanishes if and only if a closed two-form on \(M\) is exact (Theorem~\ref{thm:global-inverse}).
Finally, if \(B=0\), \(\mathcal R\) is self-adjoint and uniformly positive, and \(\nabla^{(t)}\mathcal R\) has integrable norm along a solution, then every Jacobi field along that solution is bounded (Theorem~\ref{thm:jacobi-bound}).

Section~\ref{sec:orthogonal-cartan} represents \(J^1(\mathbb R,\mathbb E^n)\) as a homogeneous space of the Galilean group and recalls Cartan geometries and the development of curves.
Section~\ref{sec:laboratory-metrics} defines Galilean manifolds and proves the existence of laboratory coordinates, in which the geodesics are the solutions of a system of second order differential equations and the connection determines the fundamental tensor.
Section~\ref{sec:local-classification} proves the fundamental theorem and defines the normal Galilean connection and its Jacobi endomorphism.
Section~\ref{sec:curvature-bianchi} gives the covariant derivatives, curvature, and Bianchi identities in laboratory coordinates, and compares the covariant derivative of a normal Galilean connection with the linear connection of Crampin, Mart\'inez, and Sarlet and with the connection of \cite{lackey2000metric}.
Section~\ref{sec:riemannian-specialization} constructs the normal Galilean connection of a second order differential equation on a Riemannian manifold and determines when its torsion coefficients \(B\), \(P\), and \(T\) vanish.
Section~\ref{sec:inverse-problems} treats the local and global inverse problems with multiplier matrix \(g\), and Section~\ref{sec:jacobi} derives the Jacobi equation and bounds its solutions.
The examples in Sections~\ref{sec:laboratory-metrics}, \ref{sec:local-classification}, \ref{sec:inverse-problems}, and~\ref{sec:jacobi} are systems of classical mechanics; in Sections~\ref{sec:inverse-problems} and~\ref{sec:jacobi} they show that the local and global inverse problems differ, and that neither the self-adjointness of \(\mathcal R\) nor the integrability of \(\|\nabla^{(t)}\mathcal R\|\) can be omitted from the bound on Jacobi fields.

\section{Background}
\label{sec:orthogonal-cartan}

We represent the first jet space as a homogeneous space of a Galilean group. In
\cite{grant2012galilean} the target \(\mathbb R^n\) is treated as an affine
space, which gives the group
\[
G_n=\left\{
\begin{pmatrix}1&0&0\\ t&1&0\\ x^i&v^i&A^i{}_j\end{pmatrix}
:\ \begin{array}{l}t\in\mathbb R,\\ x,v\in\mathbb R^n,\\ A\in GL_n\end{array}
\right\}.
\]
The jet space is
\(J^1(\mathbb R,\mathbb R^n)\cong G_n/(\mathbf 1_2\times GL_n)\),
with canonical coordinates \((t,x^i,v^i)\). The group \(G_n\) allows arbitrary linear changes of spatial frame, as the equivalence problem for systems of second order differential equations requires. The Galilean group of classical mechanics allows only orthogonal changes of spatial frame, and so has Euclidean space
\(\mathbb E^n=(\mathbb R^n,\delta)\) as its target space:
\[
\Gal_n=\left\{
\begin{pmatrix}1&0&0\\ t&1&0\\ x^i&v^i&h^i{}_j\end{pmatrix}
:\ \begin{array}{l}t\in\mathbb R,\\ x,v\in\mathbb R^n,\\ h\in O_n\end{array}
\right\}.
\]
The group acts on the jet space by
\[
(t,x,v)\stackrel{(a,b,c,h)}{\longmapsto}(a+t,\ b+ct+hx,\ c+hv),
\]
where \(a,b,c,h\) are the time, position, velocity, and orthogonal
components of a group element. The subgroup \(\mathbf 1_2\times O_n\)
is the stabilizer of the origin \((0,0,0)\). Thus
\(J^1(\mathbb R,\mathbb E^n)\cong \Gal_n/(\mathbf 1_2\times O_n)\),
again with canonical coordinates \((t,x^i,v^i)\).

The coordinates \((t,x^i,v^i)\) are time, position, and velocity, and \(J^1(\mathbb R,\mathbb E^n)\) is the setting of classical mechanics with time included. In the terminology of \cite[\S12]{souriau1997structure} it is the evolution space, as distinct from the space of motions.

The Maurer--Cartan form of \(\Gal_n\) is
\[
\omega_{\Gal_n}=
\begin{pmatrix}
0&0&0\\
dt&0&0\\
(h^{-1})^i{}_k(dx^k-v^kdt)&(h^{-1})^i{}_k\,dv^k&
(h^{-1})^i{}_k\,dh^k{}_j
\end{pmatrix}.
\]
The forms \(dx^i-v^i dt\) are called \emph{contact forms}. They encode the fact that \(v^i\)
represents the derivative of \(x^i\) with respect to \(t\): a curve
\(t\stackrel{\sigma}{\longmapsto}(t,x^i(t),v^i(t))\) has
\(\sigma^*(dx^i-v^i dt)=0\) if and only if
\(v^i(t)=\frac{dx^i}{dt}\).

The straight lines of \(J^1(\mathbb R,\mathbb E^n)\) are the integral curves of \(\frac{dx^i}{dt}=v^i\), \(\frac{dv^i}{dt}=0\), with a nondegeneracy condition:

\begin{definition}
A curve \(\sigma:I\to J^1(\mathbb R,\mathbb E^n)\) is a
\emph{straight line} if \(\sigma^*(dt)\) is nowhere zero,
\(\sigma^*(dx^i-v^i dt)=0\), and \(\sigma^*(dv^i)=0\).
\end{definition}

A \emph{Cartan geometry} \cite{sharpe1997differential} modeled on a homogeneous space
\(\varpi:G\to G/H\) of a Lie group \(G\) consists of the following data:
\begin{enumerate}
\item a smooth manifold \(M\);
\item a right principal \(H\)-bundle \(\pi:P\to M\); and,
\item a \(\mathfrak g\)-valued one-form
      \(\omega\in\Omega^1(P,\mathfrak g)\) (the ``Cartan connection'')
      which satisfies
\begin{enumerate}
\item at each \(p\in P\), \(\omega|_p:T_pP\to\mathfrak g\) is an
      isomorphism,
\item \(R_h^*\omega=\operatorname{Ad}(h^{-1})\omega\), and
\item \(\omega(Z^\dagger)=Z\) for all \(Z\in\mathfrak h\), where
      \(Z^\dagger\) is the corresponding fundamental vertical field.
\end{enumerate}
\end{enumerate}

Condition (3a), with \(H\)-equivariance, guarantees that
\(T_xM\cong\mathfrak g/\mathfrak h\) for every \(x\in M\).
For the homogeneous model \(G\to G/H\), the Maurer--Cartan form on \(G\)
satisfies the Maurer--Cartan equations. The failure of the Cartan
connection \(\omega\) to satisfy
these equations is precisely the definition of curvature:
\[
\Omega=d\omega+\frac{1}{2}[\omega,\omega].
\]
If the curvature vanishes, the Cartan geometry is locally equivalent to
\(G/H\); see \cite{sharpe1997differential}.

Each \(p\in P\) determines an isomorphism between \(T_xM\) and \(\mathfrak g/\mathfrak h\), where \(x=\pi(p)\), as follows. Write
\(VT_pP=\ker(\pi_*|_p)\) for the vertical vector space at
\(p\), so that there is a short exact sequence
\[
0\longrightarrow VT_pP\hookrightarrow T_pP
\longrightarrow T_xM\longrightarrow0.
\]
The Cartan connection provides an isomorphism
\(\omega|_p:T_pP\to\mathfrak g\), which by property (3c) restricts to
an isomorphism \(\omega|_p:VT_pP\to\mathfrak h\). Thus this short exact sequence
extends to a commutative diagram with short exact rows,
\[
\begin{array}{ccccccccc}
0&\to&VT_pP&\hookrightarrow&T_pP&\to&T_xM&\to&0\\
&&\downarrow&&\downarrow&&\downarrow&&\\
0&\to&\mathfrak h&\hookrightarrow&\mathfrak g&\to&
\mathfrak g/\mathfrak h&\to&0,
\end{array}
\]
where the first two downward maps are isomorphisms. These induce an isomorphism
\(\iota_p:T_xM\to\mathfrak g/\mathfrak h\) at each \(p\in P\).
These isomorphisms depend on \(p\), not only on \(x\). By property (3b) they are equivariant:
\(\iota_{ph}=\operatorname{Ad}(h^{-1})\iota_p\) for any
\(h\in H\), where \(\operatorname{Ad}\) here denotes the induced
action on \(\mathfrak g/\mathfrak h\). If \(T\) is a tensor on
\(\mathfrak g/\mathfrak h\) invariant under this action, for instance a
metric tensor, then \(\iota_p^*T\) is independent of the choice of
\(p\) above \(x\), and defines a tensor on \(M\).

We define geodesics through the development of curves \cite[\S3.7]{sharpe1997differential}. Let
\(\rho:(-u_0,u_0)\to P\) be a germ of a curve, and consider
\(\rho^*\omega\in\Omega^1((-u_0,u_0),\mathfrak g)\).
Then the equation
\(\tilde\rho^*\omega_G=\rho^*\omega\), where \(\omega_G\) is the
Maurer--Cartan form of \(G\), is a system of first order
ordinary differential equations for an unknown curve
\(\tilde\rho:(-u_0,u_0)\to G\). With the initial
condition \(\tilde\rho(0)=\mathbf1\), it has a unique solution,
possibly after shrinking the interval. The association
\(\rho\mapsto\tilde\rho\) is called \emph{development}.

Now, given a curve \(\sigma:(-u_0,u_0)\to M\), fix an initial frame \(p_0\in P\) above \(\sigma(0)\), and let \(\rho\) be any lift with \(\rho(0)=p_0\). The development \(\tilde\rho\) depends on the lift, but by equivariance the curve \(\tilde\sigma=\varpi\circ\tilde\rho:(-u_0,u_0)\to G/H\) depends only on \(p_0\). Replacing \(p_0\) by \(p_0h\) replaces \(\tilde\sigma\) by \(h^{-1}\tilde\sigma\), so \(\tilde\sigma\) is determined by \(\sigma\) up to the action of \(H\) on \(G/H\). We also call \(\sigma\mapsto\tilde\sigma\) development.

When \(G/H\) carries a family of straight lines invariant under \(H\), as \(J^1(\mathbb R,\mathbb E^n)\) does, a geodesic of the Cartan geometry is a curve whose development is a straight line. For the torsion-free Cartan geometry of a Riemannian metric, modeled on Euclidean space, these are the Riemannian geodesics \cite{sharpe1997differential}.

\section{Curved Galilean geometries}
\label{sec:laboratory-metrics}

In this section we define Galilean manifolds, the Cartan geometries modeled on \(J^1(\mathbb R,\mathbb E^n)\cong \Gal_n/(\mathbf 1_2\times O_n)\). Under the partial integrability condition of \cite{grant2012galilean}, a Galilean manifold has laboratory coordinates, in which its geodesics are the solutions of a system of second order differential equations and its connection determines a metric on the spatial bundle, the fundamental tensor.

\begin{definition}
A \emph{Galilean manifold} is a \((2n+1)\)-dimensional smooth manifold
\(X\), equipped with a principal \(O_n\)-bundle \(\pi:P\to X\) and a
Cartan connection with values in \(\mathfrak{gal}_n\).
\end{definition}
Here \(\mathfrak{gal}_n\) and \(\mathfrak o_n\) are the Lie algebras of
\(\Gal_n\) and \(O_n\). The Cartan connection of a Galilean manifold is called a
\emph{Galilean connection}. It has the same shape as the Maurer--Cartan form of \(\Gal_n\),
\[
\omega = \begin{pmatrix}0&0&0\\ \tau&0&0\\ \omega^a&\phi^a&\Pi^a{}_b\end{pmatrix}.
\]
The failure of the Maurer--Cartan equations defines the torsion and curvature of the connection, which for Galilean connections take the form
\[
\Omega = d\omega + \frac{1}{2}[\omega,\omega] = \begin{pmatrix}0&0&0\\ d\tau&0&0\\ T^a&\Phi^a&\Omega^a{}_b\end{pmatrix},
\]
where
\begin{align}
T^a &= d\omega^a+\Pi^a{}_b\wedge\omega^b+\phi^a\wedge\tau,
\label{eq:spatial-torsion}\\
\Phi^a &= d\phi^a+\Pi^a{}_b\wedge\phi^b,
\label{eq:velocity-torsion}\\
\Omega^a{}_b &= d\Pi^a{}_b+\Pi^a{}_c\wedge\Pi^c{}_b.
\label{eq:vertical-curvature}
\end{align}
These forms measure the failure of a Cartan geometry to be locally
equivalent to its homogeneous model. In \(J^1(\mathbb R,\mathbb E^n)\) we have canonical
coordinates \((t,x^i,v^i)\). A general Galilean manifold need not have coordinates \((t,x^i,v^i)\), even locally, together with a section \(s\) of \(P\) in which \(s^*\tau=dt\), \(s^*\omega^a=U^a{}_i(dx^i-v^i\,dt)\), and \(s^*\phi^a=U^a{}_i\,dv^i\pmod{dt,dx^j}\) for a symmetric positive definite matrix-field \(U\).
Such coordinates, called laboratory coordinates, exist locally exactly when the connection satisfies the partial integrability condition of \cite[Definition 3.2]{grant2012galilean} (Theorem~\ref{thm:laboratory} and the converse following Definition~\ref{def:laboratory}).

\begin{definition}
A Galilean connection is \emph{regular} if
\[
d\tau=0,\qquad T^a\equiv0\pmod{\omega^1,\ldots,\omega^n}.
\]
\end{definition}
The condition on \(T^a\) means that every term of this two-form contains at least one of \(\omega^1,\ldots,\omega^n\). We also call a Galilean manifold regular when its connection is regular.

\begin{theorem}[Existence of laboratory coordinates]
\label{thm:laboratory}
Let \(X\) be a regular Galilean manifold. Then about any point of \(X\)
there exists a coordinate system \((t,x^i,v^i)\) and a section \(s\)
of \(P\), such that \(s^*\tau=dt\),
\(s^*\omega^a=U^a{}_i(dx^i-v^i dt)\), and
\(s^*\phi^a=U^a{}_i\,dv^i\pmod{dt,dx^j}\), where \(U\) is a
symmetric positive definite matrix-field.
\end{theorem}

\begin{proof}
As in the proof of \cite[Theorem 3.3]{grant2012galilean}, we begin by letting
\(\tilde s\) be any local section of \(P\) about our given point.
From the regularity condition, we have
\(d(\tilde s^*\tau)=\tilde s^*d\tau=0\). By shrinking the domain
of \(\tilde s\) if necessary, we may find a function \(t\) such that
\(dt=\tilde s^*\tau\). Since \(\tau\) never vanishes, \(t\) is a
coordinate function.

Now consider the exterior differential ideal
\(\mathcal I=\langle dt,\tilde s^*\omega^a\rangle\). Since
\begin{align*}
d(\tilde s^*\omega^a)
&=\tilde s^*T^a-\tilde s^*\Pi^a{}_b\wedge\tilde s^*\omega^b
  -\tilde s^*\phi^a\wedge dt\\
&=0\pmod{\mathcal I},
\end{align*}
this ideal is integrable. As \(dt\in\mathcal I\), the integral
manifolds of \(\mathcal I\) are contained in
\(t=\mathit{constant}\) hypersurfaces. Again shrinking our neighborhood
if necessary, we may choose independent functions \(x^i\) so that
the integral manifolds of \(\mathcal I\) are given by
\(t=\mathit{constant}\) and \(x^i=\mathit{constant}\).
That is, \(\mathcal I=\langle dt,dx^i\rangle\). So,
\(\tilde s^*\omega^a\) is a combination of \(dt\) and \(dx^i\),
which we choose to write as
\[
\tilde s^*\omega^a=h^a{}_bU^b{}_i(dx^i-v^i dt)
\]
for some functions \(v^i\), using the polar decomposition of an
invertible matrix-field: \(U\) is symmetric positive definite, \(h\) is
orthogonal, and both are smooth.

Still following \cite{grant2012galilean}, we claim that
\(\{dt,dx^i,dv^i\}\) forms a local coframe, and therefore
\((t,x^i,v^i)\) is a coordinate system. First,
\[
d(\tilde s^*\omega^a)
=d(h^a{}_cU^c{}_i)\,(U^{-1})^i{}_d(h^{-1})^d{}_b\wedge\tilde s^*\omega^b
 -h^a{}_cU^c{}_i\,dv^i\wedge dt.
\]
Comparing this with~\eqref{eq:spatial-torsion} and using regularity,
we obtain
\[
(h^a{}_cU^c{}_i\,dv^i-\tilde s^*\phi^a)\wedge dt
=0\pmod{\tilde s^*\omega^b}.
\]
It follows that \(h^a{}_cU^c{}_i\,dv^i-\tilde s^*\phi^a\) is a
combination of \(dt\) and the \(\tilde s^*\omega^b\).
Yet \(dt,\tilde s^*\omega^a,\tilde s^*\phi^a\) are linearly
independent. Therefore the \(dv^i\) are linearly independent,
and independent from \(dt,dx^i\) as well.

So we have \(\tilde s^*\tau=dt\),
\(\tilde s^*\omega^a=h^a{}_cU^c{}_i(dx^i-v^i dt)\), and
\(\tilde s^*\phi^a=h^a{}_cU^c{}_i\,dv^i\pmod{dt,dx}\).
In \cite{grant2012galilean}, a general linear change of section removes
the whole coefficient matrix. Here we use the orthogonal structure group
and remove only \(h\): the section \(s=\tilde sh\) has
\(s^*\omega^a=(h^{-1})^a{}_b\tilde s^*\omega^b=U^a{}_i(dx^i-v^idt)\),
and similarly for \(s^*\phi^a\), as desired.
\end{proof}

\begin{definition}
\label{def:laboratory}
A \emph{laboratory coordinate system} is a coordinate system \((t,x^i,v^i)\) together with a section \(s\) of \(P\) such that \(s^*\tau=dt\), \(s^*\omega^a=U^a{}_i(dx^i-v^i dt)\), and \(s^*\phi^a=U^a{}_i\,dv^i\pmod{dt,dx^j}\), where \(U\) is a symmetric positive definite matrix-field.
\end{definition}

For given coordinates, the section \(s\) is unique: if \(sh\) also has these properties, then \(h^{-1}U\) is symmetric and positive definite, so \(h=1\) by the uniqueness of the polar decomposition.
The symmetry of \(U\) compares a frame index with a coordinate index, so it normalizes \(s\) relative to the coordinates and is not a tensorial condition.
In laboratory coordinates, \(s^*\phi^a-U^a{}_i\,dv^i\) is a combination of \(dt\) and the contact forms \(dx^j-v^j\,dt\), so there are unique functions \(\Gamma^i\) and \(N^i{}_j\) on the coordinate domain for which the coframe \(s^*\tau\), \(s^*\omega^a\), \(s^*\phi^a\) is
\begin{align}
s^*\tau &= dt,\nonumber\\
s^*\omega^a &= U^a{}_i(dx^i-v^i dt),\nonumber\\
s^*\phi^a &= U^a{}_i \bigl(dv^i+\Gamma^i dt+N^i{}_j(dx^j-v^jdt)\bigr). \label{eq:laboratory-forms}
\end{align}
Conversely, if \(s^*\tau=dt\), \(s^*\omega^a=U^a{}_i(dx^i-v^i\,dt)\), and \(s^*\phi^a\) has the form~\eqref{eq:laboratory-forms}, then \(d\tau=0\) and \(T^a\equiv0\pmod{\omega^b}\), since the \(dv^i\wedge dt\) terms of \(T^a\) cancel.

\begin{corollary}
\label{cor:laboratory-transitions}
Let \((t,x^i,v^i)\) and \((\bar t,\bar x^i,\bar v^i)\) be
overlapping laboratory coordinate systems associated to sections \(s,\bar s\) respectively. Locally on their overlap,
\begin{itemize}
\item \(\bar t=t+a\), for some constant \(a\);
\item \(\bar x^i=\bar x^i(t,x^j)\); and
\item \(\displaystyle
 \bar v^i=\frac{\partial\bar x^i}{\partial t}
          +\frac{\partial\bar x^i}{\partial x^j}v^j\).
\end{itemize}
When \(\bar s=sh\), we also have
\[
\bar U^a{}_i\frac{\partial\bar x^i}{\partial x^j}=(h^{-1})^a{}_bU^b{}_j.
\]
\end{corollary}

\begin{proof}
Let \(s,\bar s\) be as described. Now
\(d\bar t=\bar s^*\tau=s^*\tau=dt\), and so
\(\bar t=t+a\) locally, for some constant \(a\). Similarly,
\begin{align*}
\bar U^a{}_i(d\bar x^i-\bar v^i dt)
&=\bar s^*\omega^a\\
&=(h^{-1})^a{}_b\,s^*\omega^b\\
&=(h^{-1})^a{}_b U^b{}_j(dx^j-v^jdt).
\end{align*}
Comparing the coefficients of \(dv^j\), we find
\(\frac{\partial\bar x^i}{\partial v^j}=0\). Thus
\(\bar x^i=\bar x^i(t,x)\). Comparing the coefficients of \(dx^j\)
then gives \(\bar U^a{}_i\frac{\partial\bar x^i}{\partial x^j}=(h^{-1})^a{}_bU^b{}_j\),
and comparing the coefficients of \(dt\) gives the asserted formula
for \(\bar v^i\).
\end{proof}

In particular, changes of laboratory coordinates may depend on \(t\). Conversely, a local coordinate change of this form with invertible Jacobian matrix \(\bigl(\frac{\partial\bar x^i}{\partial x^j}\bigr)\) gives another laboratory coordinate system: one takes the polar decomposition of \(U^a{}_i\frac{\partial x^i}{\partial\bar x^j}\) to obtain the required section.

\begin{definition}
A \emph{geodesic} in a Galilean manifold is a germ of a curve
\(\sigma:(-u_0,u_0)\to X\) whose development
\(\tilde\sigma:(-u_0,u_0)\to J^1(\mathbb R,\mathbb E^n)\)
is a straight line.
\end{definition}

In laboratory coordinates the geodesics are the solutions of a system of second order differential equations. The following proposition is the analogue of \cite[Proposition 3.6]{grant2012galilean} for orthogonal frames, with the same proof.

\begin{proposition}
\label{prop:geodesic-equations}
Let \(X\) be a regular Galilean manifold and \((t,x^i,v^i)\) be
laboratory coordinates in \(X\). Write \(s^*\phi^a\) as
in~\eqref{eq:laboratory-forms}. Then the geodesics in \(X\),
parameterized by \(t\), are the solutions to the system of second
order differential equations
\begin{equation}
\frac{d^2x^i}{dt^2}
+\Gamma^i\left(t,x,\frac{dx}{dt}\right)=0.
\label{eq:geodesic-equations}
\end{equation}
\end{proposition}

\begin{proof}
By definition, the straight lines in \(J^1(\mathbb R,\mathbb E^n)\) are the curves \(\tilde\sigma\) with
\(\tilde\sigma^*(dt)\ne0\),
\(\tilde\sigma^*(dx^i-v^i dt)=0\), and
\(\tilde\sigma^*(dv^i)=0\). Up to change of basis, \(dt,dx^i-v^i dt,dv^i\) are
the \(\tau,\omega^a,\phi^a\)
components of the Maurer--Cartan form of \(\Gal_n\). By the
definition of development, a curve \(\sigma\) is a geodesic if and
only if \(\sigma^*\tau\ne0\),
\(\sigma^*(s^*\omega^a)=0\), and
\(\sigma^*(s^*\phi^a)=0\). In laboratory coordinates these equations
are \(v^i=\frac{dx^i}{dt}\) and \(\frac{dv^i}{dt}=-\Gamma^i\).
\end{proof}

\begin{corollary}
\label{cor:gamma-transitions}
If \((\bar t,\bar x,\bar v)\) is an overlapping laboratory coordinate
system whose geodesic equations are given by
\(\bar\Gamma^i(\bar t,\bar x,\bar v)\), then
\begin{equation}
\bar\Gamma^i(\bar t,\bar x,\bar v)
=\frac{\partial\bar x^i}{\partial x^j}\Gamma^j(t,x,v)
  -\frac{\partial^2\bar x^i}{\partial t^2} -2\frac{\partial^2\bar x^i}{\partial t\,\partial x^j}v^j
  -\frac{\partial^2\bar x^i}{\partial x^j\partial x^k}v^jv^k.
\label{eq:gamma-transition}
\end{equation}
\end{corollary}

\begin{proof}
Change coordinates in the geodesic equations, using
\(\bar t=t+a\) and \(\bar x^i=\bar x^i(t,x)\).
\end{proof}

Changing the section gives
\begin{align}
d\bar x^i-\bar v^i\,d\bar t
 &=\frac{\partial\bar x^i}{\partial x^j}(dx^j-v^jdt),\label{eq:x-coframe-transitions}\\
d\bar v^i+\bar\Gamma^i\,d\bar t
 +\bar N^i{}_j(d\bar x^j-\bar v^j\,d\bar t)
 &=\frac{\partial\bar x^i}{\partial x^j}\bigl(dv^j+\Gamma^jdt
                       +N^j{}_k(dx^k-v^kdt)\bigr).
\label{eq:v-coframe-transitions}
\end{align}
Indeed, differentiating
\(\bar v^i=\frac{\partial\bar x^i}{\partial t}
             +\frac{\partial\bar x^i}{\partial x^j}v^j\)
and using~\eqref{eq:gamma-transition} gives
\[
d\bar v^i+\bar\Gamma^i\,d\bar t
 =\frac{\partial\bar x^i}{\partial x^j}(dv^j+\Gamma^jdt)
   +\frac{\partial\bar v^i}{\partial x^k}(dx^k-v^kdt).
\]
Substituting into~\eqref{eq:v-coframe-transitions}
and comparing coefficients of \(dx^k-v^kdt\) gives
\begin{equation}
\bar N^i{}_j\frac{\partial\bar x^j}{\partial x^k}
=\frac{\partial\bar x^i}{\partial x^j}N^j{}_k
  -\frac{\partial\bar v^i}{\partial x^k}.
\label{eq:n-transition}
\end{equation}
Parameterizing the
geodesics by \(t\) defines a vector field \(\Xi\), with
\(\tau(\Xi)=1,\ \omega(\Xi)=\phi(\Xi)=0\). In laboratory coordinates,
\begin{equation}
\Xi=\frac{\partial}{\partial t}
 +v^i\frac{\partial}{\partial x^i}
 -\Gamma^i\frac{\partial}{\partial v^i},
\label{eq:geodesic-vector-field}
\end{equation}
and along a geodesic \(\Xi(f)=\frac{d}{dt}f(t,x(t),v(t))\) for any
function \(f\).

\begin{remark}
\label{rem:intrinsic-metric}
The action of \(O_n\) on \(\mathfrak{gal}_n/\mathfrak o_n\cong\mathbb R^{2n+1}\) preserves the inner product
\[
\tau\otimes\tau+\delta_{ab}\omega^a\otimes\omega^b+\delta_{ab}\phi^a\otimes\phi^b,
\]
so, as in Section~\ref{sec:orthogonal-cartan}, it determines a Riemannian metric \(\mathbf g\) on a Galilean manifold \(X\), for which \(s^*\tau\), \(s^*\omega^a\), \(s^*\phi^a\) is an orthonormal coframe for every local section \(s\) of \(P\).
When the connection is regular, substituting~\eqref{eq:laboratory-forms} gives, in laboratory coordinates,
\begin{equation}
\begin{split}
\mathbf g&=dt\otimes dt+
 g_{ij}(dx^i-v^idt)\otimes(dx^j-v^jdt)\\
&\quad+g_{ij}
 \bigl(dv^i+\Gamma^idt+N^i{}_k(dx^k-v^kdt)\bigr) \otimes
 \bigl(dv^j+\Gamma^jdt+N^j{}_l(dx^l-v^ldt)\bigr),
\end{split}
\label{eq:full-metric}
\end{equation}
where \(g_{ij}=\delta_{ab}U^a{}_iU^b{}_j\).
\end{remark}

The metric \(\mathbf g\) is defined on all of \(TX\). The fundamental tensor is a metric on a quotient of \(\ker\tau\). The bundle \(\ker\tau\cap\ker\omega\) is the vertical bundle, locally spanned by the \(\frac{\partial}{\partial v^i}\).

\begin{definition}
\label{def:fundamental-tensor}
The \emph{spatial bundle} of a Galilean manifold is \(\mathcal E=\ker\tau/(\ker\tau\cap\ker\omega)\). The \emph{fundamental tensor} of a Galilean manifold is the
metric on \(\mathcal E\) defined, in a local section \(s\), by
\begin{equation}
g([Y],[Z])=\delta_{ab}(s^*\omega^a)(Y)(s^*\omega^b)(Z),
\qquad Y,Z\in\ker\tau.
\label{eq:spatial-metric}
\end{equation}
In laboratory coordinates its coefficients in the frame
\([\frac{\partial}{\partial x^i}]\) are
\(g_{ij}=\delta_{ab}U^a{}_iU^b{}_j\).
\end{definition}

\begin{proposition}
The fundamental tensor is a globally well-defined, positive-definite
metric on \(\mathcal E\).
\end{proposition}
\begin{proof}
Changing a representative of \([Y]\) or \([Z]\) by a vertical vector
does not change the pairing, since the \(\omega^a\)
annihilate vertical vectors. Changing section multiplies both factors by an
orthogonal matrix and also leaves their pairing unchanged. Finally,
the kernel of \(\omega\) on \(\ker\tau\) is precisely the vertical
distribution, so
the induced metric is positive definite.
\end{proof}

On overlapping laboratory coordinate systems,
\[
g_{kl}=\bar g_{ij}
\frac{\partial\bar x^i}{\partial x^k}
\frac{\partial\bar x^j}{\partial x^l},
\]
as follows from Corollary~\ref{cor:laboratory-transitions}.

\begin{example}
\label{ex:free-fall}
In Newtonian mechanics a uniform gravitational field cannot be distinguished from an acceleration of the frame of reference; this is the equivalence principle, the starting point of Newton--Cartan geometry \cite{kunzle1972galilei,duval1985bargmann}. This example shows that the regular Galilean connection of a uniform gravitational field has this property: a change of laboratory coordinates depending on \(t\) removes the constant term from its differential equation and leaves \(\omega\), \(\phi\), and the fundamental tensor unchanged. Consider a body falling freely in a uniform gravitational field of strength \(\Gamma_0\), which obeys \(\frac{d^2x}{dt^2}+\Gamma_0=0\), where \(x\) is its height. Since \(\mathfrak o_1=0\), in the section given by \(h=1\) the forms
\[
\tau=dt,\qquad \omega=dx-v\,dt,\qquad \phi=dv+\Gamma_0\,dt,\qquad \Pi=0
\]
define a regular Galilean connection with laboratory coordinates \((t,x,v)\), \(U=1\), \(\Gamma=\Gamma_0\), and \(N=0\). The coordinates of a frame falling with the body are
\[
\bar t=t,\qquad \bar x=x+\frac{1}{2}\Gamma_0t^2,\qquad \bar v=v+\Gamma_0t,
\]
which are laboratory coordinates by Corollary~\ref{cor:laboratory-transitions}, although \(\bar x\) depends on \(t\). By~\eqref{eq:gamma-transition}, \(\bar\Gamma=\Gamma_0-\frac{\partial^2\bar x}{\partial t^2}=0\): in the falling frame every body moves on a straight line. Since
\[
d\bar x-\bar v\,d\bar t=dx-v\,dt,\qquad d\bar v+\bar\Gamma\,d\bar t=dv+\Gamma_0\,dt,
\]
the forms \(\omega\) and \(\phi\) are the same in both coordinate systems, so the metric \(\mathbf g\) and the fundamental tensor are unchanged. In contrast, \(d\bar x\otimes d\bar x=(dx+\Gamma_0t\,dt)\otimes(dx+\Gamma_0t\,dt)\) does not equal \(dx\otimes dx\).
\end{example}

\begin{example}
\label{ex:forced-oscillator}
For a system driven with period \(T_0\), the equation depends on \(t\) only modulo \(T_0\), and its solutions are naturally studied with \(t\) in the circle \(\mathbb R/T_0\mathbb Z\). This example shows that, with \(t\) in this circle, the equation defines a regular Galilean manifold whose form \(\tau\) is closed but not exact, so that it is not \(J^1(\mathbb R,M)\) with \(\tau=dt\), and that its closed geodesics are the \(T_0\)-periodic solutions. Consider a harmonic oscillator of natural frequency \(\nu_0\), driven by a periodic force of frequency \(\nu\), which obeys
\[
\frac{d^2x}{dt^2}+\nu_0^2x=\varepsilon\nu_0^2\cos\nu t,
\]
where \(\varepsilon\nu_0^2\cos\nu t\) is the driving force per unit mass \cite[\S22]{landau1976mechanics}. On \(J^1(\mathbb R,\mathbb R)\times O_1\), with \(O_1=\{\pm1\}\), the forms
\[
\tau=dt,\qquad \omega=h(dx-v\,dt),\qquad
\phi=h\bigl(dv+(\nu_0^2x-\varepsilon\nu_0^2\cos\nu t)\,dt\bigr),\qquad \Pi=0
\]
define a Galilean connection, since \(\mathfrak o_1=0\). It is regular: \(d\tau=0\) and \(T=d\omega+\phi\wedge\tau=0\). By Proposition~\ref{prop:geodesic-equations}, its geodesics parameterized by \(t\) are the solutions of the driven oscillator. Let \(T_0=\frac{2\pi}{\nu}\). The forms are invariant under the action of \(\mathbb Z\) by
\[
m\cdot(t,x,v,h)=(t+mT_0,x,v,h),
\]
which is free and properly discontinuous, so they descend to a regular Galilean connection on the quotient \(X\) of \(J^1(\mathbb R,\mathbb R)\) by \(\mathbb Z\), in which \(t\) takes values in the circle \(\mathbb R/T_0\mathbb Z\). The curve \((t,0,0)\), for \(0\le t\le T_0\), is closed in \(X\), and the integral of \(\tau\) over it is \(T_0\), so \(\tau\) is not exact.

In \(J^1(\mathbb R,\mathbb R)\) the function \(t\) increases along every geodesic, so no geodesic is closed. In \(X\), a geodesic closes after one turn around the circle of times exactly when it is a \(T_0\)-periodic solution, that is, a steady response to the driving force. Every solution is the sum of a particular solution and a free oscillation of frequency \(\nu_0\), and a free oscillation is \(T_0\)-periodic only if \(\frac{\nu_0}{\nu}\) is an integer. If \(\frac{\nu_0}{\nu}\) is not an integer, there is therefore exactly one such closed geodesic,
\[
x=\frac{\varepsilon\nu_0^2}{\nu_0^2-\nu^2}\cos\nu t.
\]
At resonance, \(\nu=\nu_0\), every solution contains the term \(\frac{\varepsilon\nu_0}{2}t\sin\nu_0t\), and no geodesic closes after one turn. If \(\nu_0=\ell\nu\) with an integer \(\ell\ge2\), the closed geodesics form a two-parameter family.
\end{example}

\section{The fundamental theorem}
\label{sec:local-classification}

In laboratory coordinates, a regular Galilean manifold determines a system of second order differential equations (Proposition~\ref{prop:geodesic-equations}) and a fundamental tensor.
We prove the local converse: given these equations and a positive-definite fundamental tensor, we classify the regular Galilean connections that realize them and determine how their torsion is constrained.
Fix coordinates \((t,x,v)\), with differential equations \(\frac{d^2x^i}{dt^2}+\Gamma^i=0\), and assume a positive-definite fundamental tensor with coefficients \(g_{ij}(t,x,v)\) is given.
Let \(U\) be a symmetric positive definite matrix-field satisfying \(g_{ij}=\delta_{ab}U^a{}_iU^b{}_j\).
Indices \(i,j,k,\ldots\) refer to coordinate frames and \(a,b,c,\ldots\) to orthonormal frames; a tensor written with indices of one type is the same tensor expressed in the corresponding frame.
Coordinate indices are lowered with \(g_{ij}\), and orthonormal indices with \(\delta_{ab}\).
Let
\[
\chi^a=U^a{}_i(dx^i-v^idt).
\]
Their exterior derivatives are
\[
d\chi^a=(dU\cdot U^{-1})^a{}_b\wedge\chi^b
       -U^a{}_i\,dv^i\wedge dt
       =\theta^a{}_b\wedge\chi^b-U^a{}_i\,dv^i\wedge dt,
\]
where we define
\[
\theta^a{}_b=(dU\cdot U^{-1})^a{}_b
             -U^a{}_i\gamma^i{}_{jk}(U^{-1})^j{}_b(dx^k-v^kdt).
\]
Here the coefficients \(\gamma^i{}_{jk}\) satisfy \(\gamma^i{}_{jk}=\gamma^i{}_{kj}\) but are otherwise arbitrary.
For an arbitrary matrix-field \(N^k{}_l\), we expand this matrix-valued one-form in the coframe \(dt\), \(dx^k-v^kdt\), \(dv^k+\Gamma^kdt+N^k{}_l(dx^l-v^ldt)\) as
\begin{equation}
\begin{split}
\theta^a{}_b
&=U^a{}_i K^i{}_j(U^{-1})^j{}_b\,dt\\
&\quad+U^a{}_i L^i{}_{jk}(U^{-1})^j{}_b(dx^k-v^kdt)\\
&\quad+U^a{}_i C^i{}_{jk}(U^{-1})^j{}_b
       \bigl(dv^k+\Gamma^kdt+N^k{}_l(dx^l-v^ldt)\bigr).
\end{split}
\label{eq:computational-decomposition}
\end{equation}
Once \(\Gamma\) and \(N\) are fixed, \(K\) and \(C\) are determined entirely by \(U\), while \(L\) also depends on the choice of \(\gamma\). Only the symmetric parts of \(K\), \(L\), and \(C\) enter the proof of the fundamental theorem, since a connection with values in \(\mathfrak o_n\) absorbs their skew parts; as in \cite[\S2]{lackey2000metric}, the following lemma gives them in terms of \(g_{ij}\).

\begin{lemma}
\label{lem:computational}
For every choice of the coefficients \(\gamma^i{}_{jk}\),
\[
K_{ij}+K_{ji}=\Xi(g_{ij}),\qquad
C_{ijk}+C_{jik}=\frac{\partial g_{ij}}{\partial v^k}.
\]
Among coefficients satisfying \(\gamma^i{}_{jk}=\gamma^i{}_{kj}\), there is a unique choice for which
\begin{equation}
L_{ijk}+L_{jik}=0.
\end{equation}
\end{lemma}

\begin{proof}
We compute the symmetric part of \(\theta^a{}_b\) by differentiating
\(g_{ij}=\delta_{ab}U^a{}_iU^b{}_j\). A direct computation gives
\begin{align*}
\theta_{ab}+\theta_{ba}
&=\biggl[\Xi(g_{ij})dt+
\left(\frac{\partial g_{ij}}{\partial x^k}
      -N^l{}_k\frac{\partial g_{ij}}{\partial v^l}
      -\gamma_{ijk}-\gamma_{jik}\right)(dx^k-v^kdt)\\
&\qquad\qquad+
\frac{\partial g_{ij}}{\partial v^k}
 \bigl(dv^k+\Gamma^kdt+N^k{}_l(dx^l-v^ldt)\bigr)\biggr]
(U^{-1})^i{}_a(U^{-1})^j{}_b.
\end{align*}
The identities for \(K\) and \(C\) follow immediately, for every \(\gamma\). Define the vector fields
\[
\frac{\delta}{\delta x^k}
=\frac{\partial}{\partial x^k}
 -N^l{}_k\frac{\partial}{\partial v^l},
\]
the analogue of the horizontal vector fields of the nonlinear connection of a Finsler metric \cite[\S2.3]{bao2012introduction}.
To obtain \(L_{ijk}+L_{jik}=0\), we must solve
\(\gamma_{ijk}+\gamma_{jik}=\frac{\delta g_{ij}}{\delta x^k}\) with \(\gamma_{ijk}\) symmetric in \(j,k\).
By Christoffel's computation, with \(\frac{\delta}{\delta x^k}\) in place of \(\frac{\partial}{\partial x^k}\), the unique solution is
\begin{equation}
\gamma_{ijk}=\frac12\left(
\frac{\delta g_{ij}}{\delta x^k}
-\frac{\delta g_{jk}}{\delta x^i}
+\frac{\delta g_{ki}}{\delta x^j}\right).
\label{eq:horizontal-christoffel}
\end{equation}
\end{proof}

\begin{theorem}[The fundamental theorem]
\label{thm:metric-general}
Consider local coordinates \((t,x,v)\), and in these coordinates
differential equations \(\frac{d^2x^i}{dt^2}+\Gamma^i=0\) and a positive-definite
fundamental tensor with coefficients \(g_{ij}\).
Let \(\mathcal E\) be the quotient of \(\ker dt\) by the span of the \(\frac{\partial}{\partial v^i}\), with the metric whose coefficients in the frame \(\bigl[\frac{\partial}{\partial x^i}\bigr]\) are \(g_{ij}\).
On the orthonormal frame bundle of \(\mathcal E\) over this coordinate neighborhood, we have the following:
\begin{enumerate}
\item There exist regular Galilean connections \((\tau, \omega^a, \phi^a, \Pi^a{}_b)\) having these
coordinates as laboratory coordinates, these equations as geodesic
equations, and \(g_{ij}\) as fundamental tensor. Their torsion forms \(T^a\) and \(\Phi^a\) are
\begin{align}
T^a&=B^a{}_b\,\tau\wedge\omega^b
 +\frac{1}{2}T^a{}_{bc}\,\omega^b\wedge\omega^c
 +M^a{}_{bc}\,\omega^b\wedge\phi^c,\label{eq:general-spatial-torsion}\\
\Phi^a&=-\mathcal R^a{}_b\,\tau\wedge\omega^b
 +(B^a{}_b+2D^a{}_b)\,\tau\wedge\phi^b\nonumber\\
&\quad+\frac{1}{2}R^a{}_{bc}\,\omega^b\wedge\omega^c
 +P^a{}_{cb}\,\omega^b\wedge\phi^c
 +\frac{1}{2}Q^a{}_{bc}\,\phi^b\wedge\phi^c,
\label{eq:general-velocity-torsion}
\end{align}
where the symmetric parts of \(B_{ab}\), \(M_{abc}\), and \(P_{abc}\) in their first two indices are determined by \(g_{ij}\), \(\Gamma^i\), and \(N^i{}_j\).
\item The tensor \(D^a{}_b\) parameterizes the possible choices
of \(\phi^a\). Once \(D^a{}_b\) is fixed, the only remaining freedom is in \(\Pi\).
\item For each fixed \(D^a{}_b\), there is precisely one such connection for every choice of the skew parts (in their first two indices) of \(B_{ab}\), \(M_{abc}\), and \(P_{abc}\).
\end{enumerate}
\end{theorem}

\begin{proof}
Let \((t,x,v)\), the differential equations, and \(g_{ij}\) be as
given. Let \(U\) be the symmetric positive definite matrix-field with \(g_{ij}=\delta_{ab}U^a{}_iU^b{}_j\).
Choose \(\gamma\) as in Lemma~\ref{lem:computational}.
The portion of the principal bundle over this coordinate system has
coordinates \((t,x,v,h)\), where \(h\in O_n\). We must define the
forms \(\tau,\omega,\phi,\Pi\) which together form a
\(\mathfrak{gal}_n\)-valued Cartan connection. Since these are to be laboratory coordinates with the prescribed
fundamental tensor, we must have \(\tau=dt\) and
\begin{equation}
\omega^a=(h^{-1})^a{}_bU^b{}_i(dx^i-v^idt).
\label{eq:reconstructed-omega}
\end{equation}
Given the differential equations, the possible choices for \(\phi\) are
\begin{equation}
\phi^a=(h^{-1})^a{}_bU^b{}_i
       \bigl(dv^i+\Gamma^idt+N^i{}_j(dx^j-v^jdt)\bigr),
\label{eq:reconstructed-phi}
\end{equation}
where \(N\) is arbitrary.

With the section \(s\) given by \(h=1\), the most general orthogonal
connection form is
\[
s^*\Pi^a{}_b=\Lambda^a{}_b\,dt+
 \Sigma^a{}_{bc}\chi^c
 +\Delta^a{}_{bc}U^c{}_i
       \bigl(dv^i+\Gamma^idt+N^i{}_j(dx^j-v^jdt)\bigr),
\]
where \(\Lambda_{ab},\Sigma_{abc},\Delta_{abc}\) are skew in their
first two indices. On the full bundle this becomes
\[
\Pi=h^{-1}(s^*\Pi)h+h^{-1}dh.
\]
The resulting forms are equivariant and reproduce the generators
of the vertical \(O_n\)-action. Their \(\tau,\omega,\phi\) components
form a coframe transverse to the fibers, so the full connection is a pointwise isomorphism.
They therefore define a Galilean connection. It is regular, because \(d\tau=0\) and the \(dv^i\wedge dt\) terms of \(T^a\) cancel.

To compute the torsion, expand the \(\mathfrak o_n\)-component \(\Pi\) in the coframe
\[
dt,\qquad dx^k-v^kdt,\qquad
dv^k+\Gamma^kdt+N^k{}_l(dx^l-v^ldt),
\]
and denote its coefficients by \(H^i{}_j\), \(H^i{}_{jk}\), and
\(V^i{}_{jk}\), respectively:
\begin{equation}
\begin{split}
\Pi^i{}_j
 &=(U^{-1})^i{}_a(s^*\Pi^a{}_b)U^b{}_j
       +(U^{-1})^i{}_a\,dU^a{}_j\\
 &=H^i{}_jdt+H^i{}_{jk}(dx^k-v^kdt)+V^i{}_{jk}\bigl(dv^k+\Gamma^kdt+N^k{}_l(dx^l-v^ldt)\bigr).
\end{split}
\label{eq:coordinate-connection}
\end{equation}
Then by the definition of \(\theta^a{}_b\),
\[
(U^{-1})^i{}_a\,dU^a{}_j
 =(U^{-1})^i{}_a\theta^a{}_bU^b{}_j+\gamma^i{}_{jk}(dx^k-v^kdt),
\]
so by the decomposition~\eqref{eq:computational-decomposition},
\begin{equation}
H^i{}_j=\Lambda^i{}_j+K^i{}_j,\qquad
H^i{}_{jk}=\Sigma^i{}_{jk}+\gamma^i{}_{jk}+L^i{}_{jk},\qquad
V^i{}_{jk}=\Delta^i{}_{jk}+C^i{}_{jk},
\label{eq:connection-identities}
\end{equation}
where \(\Lambda^i{}_j=(U^{-1})^i{}_a\Lambda^a{}_bU^b{}_j\), and
\(\Sigma^i{}_{jk},\Delta^i{}_{jk}\) are obtained from
\(\Sigma^a{}_{bc},\Delta^a{}_{bc}\) in the same way. After lowering indices,
\(\Lambda_{ij}\), \(\Sigma_{ijk}\), and \(\Delta_{ijk}\) are skew in their first two indices. Lemma~\ref{lem:computational}, together with \(\gamma_{ijk}+\gamma_{jik}=\frac{\delta g_{ij}}{\delta x^k}\), therefore gives
\begin{align}
H_{ij}+H_{ji}&=\Xi(g_{ij}),\nonumber\\*
H_{ijk}+H_{jik}&=\frac{\delta g_{ij}}{\delta x^k},
\label{eq:coordinate-metricity}\\*
V_{ijk}+V_{jik}&=\frac{\partial g_{ij}}{\partial v^k}.\nonumber
\end{align}
Equivalently, we have \(dg_{ij}=g_{ik}\Pi^k{}_j+g_{jk}\Pi^k{}_i\), which is the condition that \(s^*\Pi^a{}_b\) be skew.

As \(h,U\) only appear in pairs, put \(A=U^{-1}h\). Then
\begin{align*}
\omega^a&=(A^{-1})^a{}_i(dx^i-v^idt),\\
\phi^a&=(A^{-1})^a{}_i
       \bigl(dv^i+\Gamma^idt+N^i{}_j(dx^j-v^jdt)\bigr).
\end{align*}
Taking exterior derivatives and adding the connection terms gives
the torsion coefficients
\begin{align}
B^a{}_b&=(A^{-1})^a{}_i(H^i{}_j-N^i{}_j)A^j{}_b,\nonumber\\
T^a{}_{bc}&=(A^{-1})^a{}_i(H^i{}_{kj}-H^i{}_{jk})
                         A^j{}_bA^k{}_c,\nonumber\\
M^a{}_{bc}&=-(A^{-1})^a{}_i V^i{}_{jk}A^j{}_bA^k{}_c,\nonumber\\
D^a{}_b&=(A^{-1})^a{}_i
 \left(N^i{}_j-\frac{1}{2}\frac{\partial\Gamma^i}{\partial v^j}\right)
 A^j{}_b,\label{eq:torsion-coefficients}\\
\mathcal R^a{}_b&=(A^{-1})^a{}_i
 \left(\frac{\delta\Gamma^i}{\delta x^j}
       -\Xi(N^i{}_j)+N^i{}_lN^l{}_j\right)A^j{}_b,\nonumber\\
R^a{}_{bc}&=(A^{-1})^a{}_i
 \left(\frac{\delta N^i{}_k}{\delta x^j}
       -\frac{\delta N^i{}_j}{\delta x^k}\right)A^j{}_bA^k{}_c,\nonumber\\
P^a{}_{bc}&=(A^{-1})^a{}_i
 \left(H^i{}_{jk}-\frac{\partial N^i{}_k}{\partial v^j}\right)
 A^j{}_bA^k{}_c,\nonumber\\
Q^a{}_{bc}&=(A^{-1})^a{}_i(V^i{}_{kj}-V^i{}_{jk})
                         A^j{}_bA^k{}_c.\nonumber
\end{align}
This proves (1).

For (2), solving for \(N\) the definition of \(D^a{}_b\) in~\eqref{eq:torsion-coefficients} gives
\[
N^i{}_j=\frac{1}{2}\frac{\partial\Gamma^i}{\partial v^j}
             +A^i{}_aD^a{}_b(A^{-1})^b{}_j.
\]
Thus \(D\) determines \(N\), and hence \(\phi\), uniquely.

Finally, once \(N\) is fixed, again from~\eqref{eq:torsion-coefficients} we see \(B,M,P\) determine \(H^i{}_j\),
\(V^i{}_{jk}\), and \(H^i{}_{jk}\), respectively; their skew parts
correspond to the free choices of \(\Lambda^a{}_{b}\), \(\Delta^a{}_{bc}\), and
\(\Sigma^a{}_{bc}\). Their symmetric parts are constrained by~\eqref{eq:coordinate-metricity}:
\begin{align}
B_{ab}+B_{ba}
&=\left(\Xi(g_{ij})-g_{il}N^l{}_j-g_{jl}N^l{}_i\right)A^i{}_aA^j{}_b,
\label{eq:B-symmetric}\\
M_{abc}+M_{bac}
&=-\frac{\partial g_{ij}}{\partial v^k}
                 A^i{}_aA^j{}_bA^k{}_c,\label{eq:M-symmetric}\\
P_{abc}+P_{bac}
&=\left(\frac{\delta g_{ij}}{\delta x^k}
 -g_{il}\frac{\partial N^l{}_k}{\partial v^j}
 -g_{jl}\frac{\partial N^l{}_k}{\partial v^i}\right)
 A^i{}_aA^j{}_bA^k{}_c.\label{eq:P-symmetric}
\end{align}
Prescribing the remaining skew parts therefore determines all the
coordinate connection coefficients. Conversely, these reconstructed
coefficients satisfy metric compatibility, so
\[
s^*\Pi^a{}_b
 =U^a{}_i\Pi^i{}_j(U^{-1})^j{}_b-dU^a{}_i(U^{-1})^i{}_b
\]
is skew and defines the required connection, proving both existence and uniqueness in (3).
\end{proof}

According to this theorem, a regular Galilean connection with the given data is determined by \(D\) and the skew parts of \(B\), \(M\), and \(P\). We select one by requiring all of these to vanish. Unlike the normality convention of \cite{grant2012galilean}, this does not require \(T=0\) (see Example~\ref{ex:relativistic} below).

\begin{definition}
\label{def:normal}
A regular Galilean connection is \emph{normal} if
\[
D^a{}_b=0,\qquad B_{ab}=B_{ba},\qquad
M_{abc}=M_{bac},\qquad P_{abc}=P_{bac}.
\]
\end{definition}

These are invariant conditions. Indeed, \(T,\Phi\) vanish on vectors
tangent to the fibers of \(P\) and are orthogonally equivariant. Under a change of section \(\bar s=sh\),
their coefficient tensors satisfy, for example,
\[
\bar D_{ab}=h^d{}_ah^e{}_b D_{de},\qquad
\bar B_{ab}=h^d{}_ah^e{}_b B_{de},\qquad
\bar P_{abc}=h^d{}_ah^e{}_bh^f{}_c P_{def},
\]
and \(M\) transforms in the same way as \(P\).
Thus taking symmetric or skew parts in the first two frame indices
commutes with a change of section.

Differentiating~\eqref{eq:gamma-transition} with respect to \(v\) gives
\[
\frac{1}{2}\frac{\partial\bar\Gamma^i}{\partial\bar v^j}
             \frac{\partial\bar x^j}{\partial x^k}
 =\frac{1}{2}\frac{\partial\bar x^i}{\partial x^j}
             \frac{\partial\Gamma^j}{\partial v^k}
      -\frac{\partial\bar v^i}{\partial x^k}.
\]
Thus, by~\eqref{eq:n-transition},
\begin{equation}
\bar N^i{}_j-\frac{1}{2}\frac{\partial\bar\Gamma^i}{\partial\bar v^j}
 =\frac{\partial\bar x^i}{\partial x^k}
   \left(N^k{}_l-\frac{1}{2}\frac{\partial\Gamma^k}{\partial v^l}\right)
   \frac{\partial x^l}{\partial\bar x^j}.
\label{eq:D-invariance}
\end{equation}
Together with
\(\bar U^a{}_i\frac{\partial\bar x^i}{\partial x^j}
  =(h^{-1})^a{}_bU^b{}_j\), this proves invariance of \(D=0\)
under changes of laboratory coordinates. The remaining conditions follow from the
torsion decompositions.

\begin{corollary}[The normal Galilean connection]
\label{cor:local-normal}
For coordinates \((t,x,v)\), differential equations \(\frac{d^2x^i}{dt^2}+\Gamma^i=0\), and a positive-definite fundamental tensor with coefficients \(g_{ij}\), there is precisely one normal Galilean connection having these coordinates as laboratory coordinates, these equations as geodesic equations, and \(g_{ij}\) as fundamental tensor.
Its coordinate coefficients satisfy
\(N^i{}_j=\frac{1}{2}\frac{\partial\Gamma^i}{\partial v^j}\) and
\begin{align}
H_{ij}
&=\frac{1}{2}\left(\Xi(g_{ij})+N_{ij}-N_{ji}\right),\nonumber\\
V_{ijk}&=\frac{1}{2}\frac{\partial g_{ij}}{\partial v^k},
\label{eq:normal-coefficients}\\
H_{ijk}
&=\frac{1}{2}\left(\frac{\delta g_{ij}}{\delta x^k}
 +g_{il}\frac{\partial N^l{}_k}{\partial v^j}
 -g_{jl}\frac{\partial N^l{}_k}{\partial v^i}\right).\nonumber
\end{align}
\end{corollary}

\begin{proof}
In Theorem~\ref{thm:metric-general}, set \(D=0\) and the skew parts of \(B\), \(M\), and \(P\) to zero. The definition of \(D\) then gives \(N^i{}_j=\frac{1}{2}\frac{\partial\Gamma^i}{\partial v^j}\), and the formulas~\eqref{eq:B-symmetric}--\eqref{eq:P-symmetric} for the symmetric parts of \(B\), \(M\), and \(P\) give the stated coefficients \(H_{ij}\), \(V_{ijk}\), and \(H_{ijk}\). By~\eqref{eq:D-invariance} and the transformation laws following Definition~\ref{def:normal}, these conditions do not depend on the laboratory coordinates or the section.
\end{proof}

The coefficient \(V_{ijk}=\frac{1}{2}\frac{\partial g_{ij}}{\partial v^k}\) is analogous to the Cartan tensor \(A_{ijk}=\frac{F}{2}\frac{\partial g_{ij}}{\partial y^k}\) of a Finsler function \(F\) \cite[\S1.4]{bao2012introduction}, albeit with a different normalization: \(V_{ijk}\) corresponds to \(A_{ijk}/F\). When \(g_{ij}=\frac{\partial^2L}{\partial v^i\partial v^j}\) for a function \(L\), \(V_{ijk}\) is symmetric in all three indices.

The torsion coefficient \(\mathcal R\) of a normal Galilean connection governs the deviation of nearby geodesics (Section~\ref{sec:jacobi}).

\begin{definition}
\label{def:jacobi-endomorphism}
The \emph{Jacobi endomorphism} of a normal Galilean connection is the endomorphism \(\mathcal R\) of the spatial bundle \(\mathcal E\) whose components in the orthonormal frame of \(\mathcal E\) dual to \(s^*\omega^a\), for a local section \(s\) of \(P\), are the torsion coefficients \(\mathcal R^a{}_b\) of~\eqref{eq:general-velocity-torsion}.
\end{definition}

The torsion form \(\Phi^a\) vanishes on vectors tangent to the fibers of \(P\) and is orthogonally equivariant, so under a change of section \(\bar s=sh\) the coefficients transform as \(\bar{\mathcal R}{}_{ab}=h^d{}_ah^e{}_b\mathcal R_{de}\), and \(\mathcal R\) is well defined. Substituting \(N^i{}_j=\frac{1}{2}\frac{\partial\Gamma^i}{\partial v^j}\) into~\eqref{eq:torsion-coefficients}, and using \(\frac{\delta\Gamma^i}{\delta x^j}=\frac{\partial\Gamma^i}{\partial x^j}-2N^i{}_lN^l{}_j\), gives
\begin{equation}
\mathcal R^a{}_b
 =(A^{-1})^a{}_i
   \left(\frac{\partial\Gamma^i}{\partial x^j}
     -\frac{1}{4}\frac{\partial\Gamma^i}{\partial v^k}
                   \frac{\partial\Gamma^k}{\partial v^j}
     -\frac{1}{2}\Xi\left(\frac{\partial\Gamma^i}{\partial v^j}\right)
   \right)A^j{}_b.
\label{eq:jacobi-endomorphism}
\end{equation}
Its coordinate components \(\mathcal R^i{}_j\) depend only on the differential equations, not on the fundamental tensor, and are the components of the Jacobi endomorphism of a system of second order differential equations \cite{crampin1996linear}.

\begin{example}
\label{ex:relativistic}
In special relativity the inertia of a particle grows with its speed, so the second derivative of its Lagrangian with respect to \(v\) depends on \(v\). Taking this second derivative as fundamental tensor gives a regular Galilean manifold whose fundamental tensor depends on \(v\). The orthogonal analogue of the normality convention of \cite{grant2012galilean}, which requires \(T=0\), excludes every fundamental tensor that depends on \(v\), including this one. This example shows that normality as defined here admits it. A free particle of rest mass \(m\) in special relativity, described in coordinate time \(t\) rather than proper time, has the Lagrangian \(-mc^2\sqrt{1-\frac{v^2}{c^2}}\) \cite[\S8]{landau1975fields}, where \(c\) is the speed of light. Divided by \(m\), its Euler--Lagrange equation is \(\frac{d^2x}{dt^2}=0\), and
\[
g_{11}=\frac{\partial^2}{\partial v^2}\Bigl(-c^2\sqrt{1-\frac{v^2}{c^2}}\Bigr)=\Bigl(1-\frac{v^2}{c^2}\Bigr)^{-3/2},
\]
the longitudinal mass of the particle per unit rest mass. Let \(X\) be the open subset \(\lvert v\rvert<c\) of \(J^1(\mathbb R,\mathbb R)\), with coordinates \((t,x,v)\), and prescribe the equation \(\frac{d^2x}{dt^2}=0\) and the fundamental tensor \(g_{11}\). The orthonormal frame bundle of \(\mathcal E\) is trivial; use the section \(s\) given by \(h=1\), and write \(\omega,\phi,\Pi\) for the pullbacks of the Cartan connection forms by \(s\). The normal Galilean connection has
\[
N=0,\qquad \Pi=0,\qquad
\omega=g_{11}^{1/2}(dx-v\,dt),\qquad
\phi=g_{11}^{1/2}\,dv,
\]
since \(N=\frac{1}{2}\frac{\partial\Gamma}{\partial v}=0\) and \(\mathfrak o_1=0\). Since \(\Pi=0\), by the definition of \(T\),
\[
T=d\omega+\phi\wedge\tau=\frac{1}{2}g_{11}^{-1/2}\frac{\partial g_{11}}{\partial v}\,dv\wedge(dx-v\,dt)
 =-\frac{1}{2}g_{11}^{-3/2}\frac{\partial g_{11}}{\partial v}\,\omega\wedge\phi.
\]
Since \(d\phi=\frac{1}{2}g_{11}^{-1/2}\frac{\partial g_{11}}{\partial v}\,dv\wedge dv=0\), we have \(\Phi=d\phi+\Pi\wedge\phi=0\), so \(D=P=0\), and \(T\) has no \(\tau\wedge\omega\) term, so \(B=0\). The coefficient
\[
M_{111}=-\frac{1}{2}g_{11}^{-3/2}\frac{\partial g_{11}}{\partial v}=-\frac{3v}{2c^2}\Bigl(1-\frac{v^2}{c^2}\Bigr)^{-1/4}
\]
is nonzero for \(v\ne0\). The connection is normal, but \(T\ne0\), so the orthogonal analogue of the normality convention \(T=0\) of \cite{grant2012galilean} excludes relativistic free motion.
\end{example}

\section{Covariant derivatives, curvature, and Bianchi identities}
\label{sec:curvature-bianchi}

In this section we write the covariant derivative and the curvature of a regular Galilean connection in laboratory coordinates, compare the covariant derivative of a normal Galilean connection with those of \cite{crampin1996linear} and \cite{lackey2000metric}, and record the Bianchi identities.
Sections~\ref{sec:riemannian-specialization}--\ref{sec:jacobi} use the temporal covariant derivative \(\nabla^{(t)}\), this comparison, the curvature components \(R^i{}_{jkl}\), and the Jacobi endomorphism \(\mathcal R\); in Sections~\ref{sec:inverse-problems} and~\ref{sec:jacobi}, \(R^i{}_{jkl}\) is the Riemann curvature tensor of the metric.
The other components are recorded for comparison with Finsler geometry and with the geometry of time-dependent second order differential equations.

The adjoint action of \(O_n\subset \Gal_n\) on
\(\mathfrak{gal}_n/\mathfrak{o}_n\) is
\(h\cdot(t,x,v)=(t,hx,hv)\).
Since this action fixes the \(t\)-component, the \(\mathfrak o_n\)-component \(\Pi\) acts trivially on \(\tau\), and the component \(\tau(Y)\) of a vector field \(Y\) on \(X\) is a function, differentiated without a connection term.
On \(X\), we call vectors in the span of the
\(\frac{\delta}{\delta x^i}\) horizontal and those in the span of the
\(\frac{\partial}{\partial v^i}\) vertical. A horizontal or vertical vector field \(Y\) on \(X\) has components \(\xi^a\) on \(P\): for any vector field \(\tilde Y\) on \(P\) projecting to \(Y\), \(\xi^a=\omega^a(\tilde Y)\) if \(Y\) is horizontal and \(\xi^a=\phi^a(\tilde Y)\) if \(Y\) is vertical. Since \(\omega^a,\phi^a\) vanish on vectors tangent to the fibers, \(\xi^a\) does not depend on the choice of \(\tilde Y\). In both cases the covariant derivative is
\[
(\nabla\xi)^a=d\xi^a+\Pi^a{}_b\xi^b.
\]
If \(Y=\xi^i(t,x,v)\frac{\delta}{\delta x^i}\) or \(Y=\xi^i(t,x,v)\frac{\partial}{\partial v^i}\) in laboratory coordinates, then in both cases \(\xi^a=(A^{-1})^a{}_i\xi^i\). Using the coordinate coefficients of the \(\mathfrak{o}_n\)-component \(\Pi\)
in~\eqref{eq:coordinate-connection}, we obtain
\begin{align}
(\nabla^{(t)}\xi)^a
 &= (A^{-1})^a{}_i\left(\Xi(\xi^i)+N^i{}_j\xi^j\right)
       +B^a{}_b\xi^b,\label{eq:nabla-temporal}\\
(\nabla^{(h)}_c\xi)^a
 &= (A^{-1})^a{}_i
 \left(\frac{\delta\xi^i}{\delta x^k}
       +\frac{\partial N^i{}_k}{\partial v^j}\xi^j\right)A^k{}_c
       +P^a{}_{bc}\xi^b,\label{eq:nabla-horizontal}\\
(\nabla^{(v)}_c\xi)^a
 &= (A^{-1})^a{}_i
       \frac{\partial\xi^i}{\partial v^k}A^k{}_c
       -M^a{}_{bc}\xi^b.\label{eq:nabla-vertical}
\end{align}
Here \(\nabla^{(t)}\xi\), \(\nabla^{(h)}_c\xi\), and \(\nabla^{(v)}_c\xi\) are \(\nabla\xi\) evaluated on lifts to \(P\) of \(\Xi\), \(A^k{}_c\frac{\delta}{\delta x^k}\), and \(A^k{}_c\frac{\partial}{\partial v^k}\); since \(\nabla\xi\) vanishes on vectors tangent to the fibers of \(P\), they do not depend on the lifts.
They are the temporal, horizontal, and vertical components of the covariant derivative.
The last terms of~\eqref{eq:nabla-temporal}--\eqref{eq:nabla-vertical} contain the torsion coefficients \(B^a{}_b\), \(P^a{}_{bc}\), and \(-M^a{}_{bc}\).

When \(D=0\), the remaining terms are the covariant derivative of the linear connection \(\mathring\nabla\) that Crampin, Mart\'inez, and Sarlet attach to the differential equations alone \cite{crampin1996linear}, which is essentially that of Massa and Pagani \cite{massa1994jet}.
In the directions \(\Xi\), \(\frac{\delta}{\delta x^k}\), and \(\frac{\partial}{\partial v^k}\), its coefficients are \(N^i{}_j\), \(\frac{\partial N^i{}_k}{\partial v^j}\), and \(0\) \cite[p.~232]{crampin1996linear}, \cite[(2.39)]{massa1994jet}, and \(\mathring\nabla_\Xi\) is the dynamical covariant derivative \cite{sarlet1995derivations}, \cite[(5)]{bucataru2007metric}.
The covariant derivatives of the fundamental tensor with respect to \(\mathring\nabla\) are
\begin{equation}
\begin{split}
(\mathring\nabla_\Xi g)_{ij}&=\Xi(g_{ij})-g_{il}N^l{}_j-g_{jl}N^l{}_i,\\
(\mathring\nabla_kg)_{ij}&=\frac{\delta g_{ij}}{\delta x^k}-g_{il}\frac{\partial N^l{}_k}{\partial v^j}-g_{jl}\frac{\partial N^l{}_k}{\partial v^i},\\
(\mathring\nabla^{(v)}_kg)_{ij}&=\frac{\partial g_{ij}}{\partial v^k}.
\end{split}
\label{eq:berwald-metric-derivatives}
\end{equation}
By~\eqref{eq:B-symmetric}--\eqref{eq:P-symmetric}, the symmetric parts of \(B\), \(P\), and \(-M\) are these covariant derivatives of \(g\), expressed in the orthonormal frame of \(\mathcal E\) dual to \(s^*\omega^a\).
For a normal Galilean connection, \(B\), \(M\), and \(P\) are symmetric in their first two indices, so \(B\), \(P\), and \(-M\) are half of~\eqref{eq:berwald-metric-derivatives}.
Thus the covariant derivative of the normal Galilean connection is the unique covariant derivative that preserves \(g\) and differs from \(\mathring\nabla\) by terms self-adjoint with respect to \(g\): if \(\mathring\nabla+S\) preserves \(g\), then \(\mathring\nabla g=g(S\,\cdot,\cdot)+g(\cdot,S\,\cdot)\), which determines \(S=\frac12g^{-1}\mathring\nabla g\) when \(S\) is self-adjoint.
In particular, \(\nabla^{(t)}-B\) is the dynamical covariant derivative and \(B\) is half the dynamical covariant derivative of \(g\).

The connection obtained by Cartan's method of equivalence in \cite[\S2]{lackey2000metric} is the regular Galilean connection with \(B=0\), \(T=0\), \(D_{ab}=D_{ba}\), and \(M_{abc}=M_{bac}\).
Indeed, by~\eqref{eq:torsion-coefficients} and~\eqref{eq:B-symmetric}, \(B=0\) and \(D_{ab}=D_{ba}\) hold exactly when \(H^i{}_j=N^i{}_j\) and
\begin{equation}
N^i{}_j=\frac{1}{2}\frac{\partial\Gamma^i}{\partial v^j}+\frac{1}{2}g^{ik}(\mathring\nabla_\Xi g)_{kj},
\label{eq:metric-nonlinear-connection}
\end{equation}
the metric nonlinear connection of \(g\), obtained in \cite[\S2]{lackey2000metric} and, for \(g\) independent of \(t\), in \cite[Theorem~2.2]{bucataru2007metric}.
The condition \(T=0\) holds exactly when \(H_{ijk}=H_{ikj}\), which by~\eqref{eq:coordinate-metricity} and Christoffel's computation gives \(H_{ijk}=\gamma_{ijk}\) as in~\eqref{eq:horizontal-christoffel}.
By~\eqref{eq:normal-coefficients}, \(g_{ik}N^k{}_j\) for the nonlinear connection~\eqref{eq:metric-nonlinear-connection} equals the coefficient \(H_{ij}\) of the normal Galilean connection, so the two connections have the same temporal covariant derivative \(\nabla^{(t)}\).
Table~\ref{tab:comparison} lists the objects of the normal Galilean connection that appear in the literature.

\begin{table}[t]
\centering
\small
\begin{tabular}{@{}p{0.3\textwidth}p{0.62\textwidth}@{}}
\toprule
Normal Galilean connection & In the literature\\
\midrule
\(N^i{}_j=\frac12\frac{\partial\Gamma^i}{\partial v^j}\) & nonlinear connection of the differential equations \cite[p.~229]{crampin1996linear}, \cite[(2.40)]{massa1994jet}\\
first terms of~\eqref{eq:nabla-temporal}--\eqref{eq:nabla-vertical} & linear connection \(\mathring\nabla\) \cite[p.~232]{crampin1996linear}, \cite[(2.39)]{massa1994jet}\\
\(\nabla^{(t)}-B=\mathring\nabla_\Xi\) & dynamical covariant derivative \cite{sarlet1995derivations}, \cite[(5)]{bucataru2007metric}; covariant derivative of the theory of Kosambi, Cartan, and Chern \cite[(1.7)]{antonelli2001invariants}\\
\(\nabla^{(t)}\) & dynamical covariant derivative of the metric nonlinear connection~\eqref{eq:metric-nonlinear-connection} \cite[Theorem~2.2]{bucataru2007metric} (for \(g\) independent of \(t\)); temporal covariant derivative of the connection of \cite[\S2]{lackey2000metric}\\
\(\mathcal R\) & Jacobi endomorphism \cite[p.~231]{crampin1996linear}, \cite[(112)]{sarlet1995derivations}; second invariant of the theory of Kosambi, Cartan, and Chern \cite[(3.4)]{antonelli2001invariants}, \cite{bucataru2010helmholtz}\\
\(B=\frac12g^{-1}\mathring\nabla_\Xi g\) & \(B=0\) is the first Helmholtz condition \cite[(4.7)]{douglas1941solution} (Section~\ref{sec:inverse-problems})\\
\bottomrule
\end{tabular}
\caption{Objects of the normal Galilean connection and their counterparts in the literature.}
\label{tab:comparison}
\end{table}

The \(\mathfrak{o}_n\)-component \(\Omega^a{}_b\) of the curvature has five components, defined by
\begin{align}
\Omega^a{}_b
 &=I^a{}_{bc}\,\tau\wedge\omega^c
   +J^a{}_{bc}\,\tau\wedge\phi^c\nonumber\\
 &\quad+\frac{1}{2}R^a{}_{bcd}\,\omega^c\wedge\omega^d
   +P^a{}_{bcd}\,\omega^c\wedge\phi^d
   +\frac{1}{2}Q^a{}_{bcd}\,\phi^c\wedge\phi^d.
\label{eq:curvature-components}
\end{align}
We compute these components using the coordinate expression for \(\Pi\).
With \(\mathcal R^i{}_j,R^i{}_{jk}\) the coordinate components of the torsion coefficients \(\mathcal R^a{}_b,R^a{}_{bc}\), the vector fields
\(\Xi,\frac{\delta}{\delta x^i},\frac{\partial}{\partial v^i}\) satisfy
\begin{align*}
\left[\Xi,\frac{\delta}{\delta x^j}\right]
 &=N^i{}_j\frac{\delta}{\delta x^i}+\mathcal R^i{}_j\frac{\partial}{\partial v^i},\\
\left[\Xi,\frac{\partial}{\partial v^j}\right]
 &=-\frac{\delta}{\delta x^j}+
   \left(\frac{\partial\Gamma^i}{\partial v^j}-N^i{}_j\right)
     \frac{\partial}{\partial v^i},\\
\left[\frac{\delta}{\delta x^j},\frac{\delta}{\delta x^k}\right]&=-R^i{}_{jk}\frac{\partial}{\partial v^i},\\
\left[\frac{\delta}{\delta x^j},\frac{\partial}{\partial v^k}\right]
 &=\frac{\partial N^i{}_j}{\partial v^k}\frac{\partial}{\partial v^i},\\
 \left[\frac{\partial}{\partial v^j},\frac{\partial}{\partial v^k}\right] &= 0.
\end{align*}
In the coordinate frame,
\(\Omega^i{}_j=d\Pi^i{}_j+\Pi^i{}_k\wedge\Pi^k{}_j\).
Decompose it as in~\eqref{eq:curvature-components}, using
\[
dt,\qquad dx^i-v^idt,\qquad
dv^i+\Gamma^idt+N^i{}_j(dx^j-v^jdt),
\]
with coefficients
\(I^i{}_{jk},J^i{}_{jk},R^i{}_{jkl},
P^i{}_{jkl},Q^i{}_{jkl}\).
Evaluation on these vector fields gives
\begin{align}
I^i{}_{jk}
 &=\Xi(H^i{}_{jk})-\frac{\delta H^i{}_j}{\delta x^k}
   +H^i{}_lH^l{}_{jk}-H^i{}_{lk}H^l{}_j
   -H^i{}_{jl}N^l{}_k-V^i{}_{jl}\mathcal R^l{}_k,\nonumber\\
J^i{}_{jk}
 &=\Xi(V^i{}_{jk})-\frac{\partial H^i{}_j}{\partial v^k}
   +H^i{}_lV^l{}_{jk}-V^i{}_{lk}H^l{}_j+H^i{}_{jk}\nonumber\\
 &\quad-V^i{}_{jl}
       \left(\frac{\partial\Gamma^l}{\partial v^k}-N^l{}_k\right),\nonumber\\
R^i{}_{jkl}
 &=\frac{\delta H^i{}_{jl}}{\delta x^k}
   -\frac{\delta H^i{}_{jk}}{\delta x^l}
   +H^i{}_{mk}H^m{}_{jl}-H^i{}_{ml}H^m{}_{jk}
   +V^i{}_{jm}R^m{}_{kl},
\label{eq:adapted-curvature}\\
P^i{}_{jkl}
 &=\frac{\delta V^i{}_{jl}}{\delta x^k}
   -\frac{\partial H^i{}_{jk}}{\partial v^l}
   +H^i{}_{mk}V^m{}_{jl}-V^i{}_{ml}H^m{}_{jk}
   -V^i{}_{jm}\frac{\partial N^m{}_k}{\partial v^l},\nonumber\\
Q^i{}_{jkl}
 &=\frac{\partial V^i{}_{jl}}{\partial v^k}
   -\frac{\partial V^i{}_{jk}}{\partial v^l}
   +V^i{}_{mk}V^m{}_{jl}-V^i{}_{ml}V^m{}_{jk}.\nonumber
\end{align}
The first four terms of \(R^i{}_{jkl}\) have the form of the curvature of a linear connection with coefficients \(H^i{}_{jk}\), with \(\frac{\delta}{\delta x^k}\) in place of \(\frac{\partial}{\partial x^k}\); the last term is \(V^i{}_{jm}R^m{}_{kl}\). The components
\(P^a{}_{bcd}\) and \(Q^a{}_{bcd}\) are analogous to the curvatures \(P\) and \(Q\) of a Finsler connection \cite[\S3.1]{bao2012introduction}.
The two components \(I^a{}_{bc},J^a{}_{bc}\) are the coefficients of the terms containing \(\tau\).
For a Lagrangian depending on \(t\), they are parts of the curvature forms in the structure equations of \cite[\S5.8]{miron2003geometry}; nonlinear connections depending on \(t\) are treated in \cite[Chapter~6]{antonelli2003kcc}.

Since \(d\tau=0\), exterior differentiation of the definitions of \(T^a\), \(\Phi^a\), and \(\Omega^a{}_b\) gives the Bianchi identities
\begin{align}
dT^a+\Pi^a{}_b\wedge T^b
 &=\Omega^a{}_b\wedge\omega^b+\Phi^a\wedge\tau,\label{eq:bianchi-spatial-torsion}\\
d\Phi^a+\Pi^a{}_b\wedge\Phi^b
 &=\Omega^a{}_b\wedge\phi^b,\label{eq:bianchi-velocity-torsion}\\
d\Omega^a{}_b+\Pi^a{}_c\wedge\Omega^c{}_b
                  -\Omega^a{}_c\wedge\Pi^c{}_b&=0.\label{eq:bianchi-curvature}
\end{align}
Each is an equality of three-forms, whose coefficients give tensor
identities. For example, by the definitions of \(M^a{}_{bc}\) and \(Q^a{}_{bc}\) in~\eqref{eq:torsion-coefficients},
\[
M^a{}_{bc}-M^a{}_{cb}=Q^a{}_{bc}.
\]
The other component identities are not used in this paper; for their analogue in Finsler spaces, see \cite[\S2.6]{miron2001geometry} or \cite[\S\S3.2--3.5]{bao2012introduction}.

\section{Differential equations on Riemannian manifolds}
\label{sec:riemannian-specialization}

In this section we take the fundamental tensor to be a Riemannian metric \(g\) on a manifold \(M\), independent of \(t\) and \(v\). We show that every second order differential equation on \(M\) parameterized by \(t\) is the vector field of a unique normal Galilean connection on the jet space \(J^1(\mathbb R,M)\) (Proposition~\ref{prop:global-normal}), and we determine when its torsion coefficients \(B\), \(P\), and \(T\) vanish (Proposition~\ref{prop:riemannian-quadratic}). The coefficients \(N^i{}_j=\frac{1}{2}\frac{\partial\Gamma^i}{\partial v^j}\) depend only on the differential equation; they form the canonical nonlinear connection of the geometric theory of second order differential equations \cite{crampin1996linear,bucataru2011geometric}.
The metric \(g\) extends this nonlinear connection to a Galilean connection, and \(B=0\) holds exactly when \(\Gamma^i=\gamma^i{}_{jk}v^jv^k+W^i{}_jv^j+F^i\), where \(\gamma^i{}_{jk}\) are the Levi-Civita coefficients of \(g\) and \(g_{il}W^l{}_j\) is skew in \(i,j\); this condition is used throughout Sections~\ref{sec:inverse-problems} and~\ref{sec:jacobi}.

Let \((M,g)\) be a smooth paracompact Riemannian manifold and \(X=J^1(\mathbb R,M)\). A coordinate system \((x^i)\) on an open set \(\mathcal O\subset M\) induces \emph{adapted coordinates} \((t,x^i,v^i)\) on the part of \(X\) over \(\mathcal O\); in them the projections \(p:X\to\mathbb R\times M\) and \(q=\operatorname{pr}_M\circ p:X\to M\) are \(p(t,x,v)=(t,x)\) and \(q(t,x,v)=x\), and the contact forms are \(dx^i-v^i\,dt\). Write \(\tau=p^*dt\). A change of coordinates on \(M\) transforms the contact forms by the Jacobian matrix \(\frac{\partial\bar x^i}{\partial x^j}\), so their span generates a globally defined ideal \(\mathcal C\), the contact ideal.

\begin{definition}
\label{def:second-order-equation}
A \emph{second order differential equation} on \(M\) is a vector field \(\mathscr S\) on \(X\) such that \(\mathcal C(\mathscr S)=0\) and \(\tau(\mathscr S)\ne0\). It is \emph{parameterized by \(t\)} if \(\tau(\mathscr S)=1\).
\end{definition}

A nowhere-vanishing one-form on \(\mathbb R\) is \(f(t)\,dt\) with \(f\) never zero, so the condition \(\tau(\mathscr S)\ne0\) does not depend on the choice of \(dt\). Multiplying by \(\frac{1}{\tau(\mathscr S)}\) parameterizes any second order differential equation by \(t\).

For \(Y\in\ker\tau\), the values \((dx^i-v^i\,dt)(Y)\) are the components of \(q_*Y\), whose kernel is the vertical bundle spanned by the \(\frac{\partial}{\partial v^i}\); thus \(q_*\) identifies the spatial bundle \(\mathcal E\) with \(q^*TM\), and we give \(\mathcal E\) the metric \(q^*g\), with coefficients \(g_{ij}(x)\). Let \(P=q^*F(M)\) be the pullback of the orthonormal frame bundle of \(g\), with projection \(\pi:P\to X\). A point of \(P\) is \((t,x,v,e)\), with \(e\) an orthonormal frame of \(T_xM\). In adapted coordinates its components are \(A=U^{-1}h\) with \(h\in O_n\), so \(P\) has coordinates \((t,x,v,h)\). With \(\chi^a=U^a{}_i(dx^i-v^i\,dt)\), set
\[
\omega^a=(h^{-1})^a{}_b\chi^b=(A^{-1})^a{}_i(dx^i-v^i\,dt).
\]
A change of coordinates multiplies \(A\) and the contact forms by the same Jacobian matrix, so \(\omega\) is defined on all of \(P\).

Let \(\mathcal O\subset M\) be a coordinate domain. In adapted coordinates, a second order differential equation on \(q^{-1}(\mathcal O)\) parameterized by \(t\) has the form
\[
\mathscr S_{\mathcal O}=\frac{\partial}{\partial t}+v^i\frac{\partial}{\partial x^i}-\Gamma^i(t,x,v)\frac{\partial}{\partial v^i},
\qquad
\frac{d^2x^i}{dt^2}+\Gamma^i\left(t,x,\frac{dx}{dt}\right)=0.
\]
Here \(\Gamma^i\) may depend on \(t,x,v\), while \(g_{ij}\) depends only on \(x\). By Corollary~\ref{cor:local-normal}, there is precisely one normal Galilean connection having the adapted coordinates as laboratory coordinates, these equations as geodesic equations, and \(g\) as fundamental tensor; its vector field \(\Xi\) is \(\mathscr S_{\mathcal O}\).

For this connection \(D=0\), so
\(N^i{}_j=\frac{1}{2}\frac{\partial\Gamma^i}{\partial v^j}\).
Substituting \(g_{ij}=g_{ij}(x)\) into~\eqref{eq:normal-coefficients},
the coefficients of the \(\mathfrak{o}_n\)-component \(\Pi\) are
\begin{align}
H_{ij}&=\frac{1}{2}\left(v^k\frac{\partial g_{ij}}{\partial x^k}+N_{ij}-N_{ji}\right),\\
V_{ijk}&=0,\\
H_{ijk}&=\frac{1}{2}\left(\frac{\partial g_{ij}}{\partial x^k}
 +g_{il}\frac{\partial N^l{}_k}{\partial v^j}
 -g_{jl}\frac{\partial N^l{}_k}{\partial v^i}\right).
\label{eq:riemannian-connection-coefficients}
\end{align}
By normality, \(B_{ab}\), \(M_{abc}\), and \(P_{abc}\) are symmetric in their first two indices, and
\begin{align}
B_{ab}+B_{ba}
 &=\left(v^k\frac{\partial g_{ij}}{\partial x^k}
    -\frac{1}{2}g_{il}\frac{\partial\Gamma^l}{\partial v^j}
    -\frac{1}{2}g_{jl}\frac{\partial\Gamma^l}{\partial v^i}\right)
        A^i{}_aA^j{}_b,\label{eq:riemannian-B-coefficients}\\
M_{abc}+M_{bac}&=0,\nonumber\\
P_{abc}+P_{bac}
 &=\left(\frac{\partial g_{ij}}{\partial x^k}
    -\frac{1}{2}g_{il}
       \frac{\partial^2\Gamma^l}{\partial v^j\partial v^k}
    -\frac{1}{2}g_{jl}
       \frac{\partial^2\Gamma^l}{\partial v^i\partial v^k}\right)
        A^i{}_aA^j{}_bA^k{}_c.\label{eq:riemannian-P-coefficients}
\end{align}

\begin{proposition}
\label{prop:riemannian-quadratic}
For the normal Galilean connection determined by \(g_{ij}(x)\) and \(\mathscr S_{\mathcal O}\):
\begin{enumerate}
\item \(B^a{}_b=0\) implies \(P^a{}_{bc}=0\), and \(P^a{}_{bc}=0\) implies \(T^a{}_{bc}=0\).
\item \(T^a{}_{bc}=0\) if and only if
\begin{equation}
P^a{}_{bc}=(A^{-1})^a{}_i
 \left(\gamma^i{}_{jk}
       -\frac{1}{2}\frac{\partial^2\Gamma^i}
                         {\partial v^j\partial v^k}\right)
 A^j{}_bA^k{}_c,
\label{eq:riemannian-torsion-free-P}
\end{equation}
where \(\gamma^i{}_{jk}\) are the Levi-Civita coefficients of \(g\).
\item \(P^a{}_{bc}=0\) if and only if there are functions \(W^i{}_j(t,x)\) and \(F^i(t,x)\) such that
\begin{equation}
\Gamma^i=\gamma^i{}_{jk}v^jv^k+W^i{}_j(t,x)v^j+F^i(t,x).
\label{eq:quadratic-equation}
\end{equation}
\item If \(\Gamma^i\) has the form~\eqref{eq:quadratic-equation}, then \(B=0\) if and only if \(g_{il}W^l{}_j+g_{jl}W^l{}_i=0\).
\end{enumerate}
\end{proposition}

\begin{proof}
Since \(N^i{}_j=\frac{1}{2}\frac{\partial\Gamma^i}{\partial v^j}\), the torsion coefficients~\eqref{eq:torsion-coefficients} satisfy
\begin{align*}
T^a{}_{bc}&=(A^{-1})^a{}_i(H^i{}_{kj}-H^i{}_{jk})A^j{}_bA^k{}_c=P^a{}_{cb}-P^a{}_{bc},\\
P^a{}_{bc}&=(A^{-1})^a{}_i\left(H^i{}_{jk}-\frac{1}{2}\frac{\partial^2\Gamma^i}{\partial v^j\partial v^k}\right)A^j{}_bA^k{}_c,
\end{align*}
and by normality \(B_{ab}=B_{ba}\) and \(P_{abc}=P_{bac}\).

Suppose \(B=0\). Then the expression in parentheses in~\eqref{eq:riemannian-B-coefficients} vanishes for every \(v\), and so does its derivative with respect to \(v^k\), which is the expression in parentheses in~\eqref{eq:riemannian-P-coefficients}. Hence \(P_{abc}+P_{bac}=0\), and since \(P_{abc}=P_{bac}\), \(P^a{}_{bc}=0\). If \(P^a{}_{bc}=0\), then \(T^a{}_{bc}=P^a{}_{cb}-P^a{}_{bc}=0\). This proves part (1).

Suppose next that \(T^a{}_{bc}=0\), so that \(H^i{}_{jk}=H^i{}_{kj}\). Since \(g\) depends only on \(x\), the metricity relation~\eqref{eq:coordinate-metricity} is
\[
g_{il}H^l{}_{jk}+g_{jl}H^l{}_{ik}
 =\frac{\partial g_{ij}}{\partial x^k},
\]
so by Christoffel's computation \(H^i{}_{jk}=\gamma^i{}_{jk}\). Substituting this into the formula for \(P^a{}_{bc}\) gives~\eqref{eq:riemannian-torsion-free-P}. Conversely, if \(P^a{}_{bc}\) has this form, the formula for \(P^a{}_{bc}\) gives \(H^i{}_{jk}=\gamma^i{}_{jk}\), which is symmetric in \(j,k\), so \(T^a{}_{bc}=0\). This proves part (2).

If \(P^a{}_{bc}=0\), then \(T^a{}_{bc}=0\) by part (1), and by part (2)
\[
\frac{\partial^2\Gamma^i}{\partial v^j\partial v^k}
 =2\gamma^i{}_{jk},
\]
so \(\Gamma^i-\gamma^i{}_{jk}v^jv^k\) has zero second derivative in \(v\). Integrating twice in \(v\) gives~\eqref{eq:quadratic-equation}. Conversely, if \(\Gamma^i\) has this form, then
\(N^i{}_j=\gamma^i{}_{jk}v^k+\frac12W^i{}_j\).
Substitution into~\eqref{eq:riemannian-connection-coefficients}, using
metric compatibility of \(\gamma\), gives \(H^i{}_{jk}=\gamma^i{}_{jk}\), and then the formula for \(P^a{}_{bc}\) gives \(P^a{}_{bc}=0\). This proves part (3).

If \(\Gamma^i=\gamma^i{}_{jk}v^jv^k+W^i{}_jv^j+F^i\), substitution into~\eqref{eq:riemannian-B-coefficients} gives
\[
B_{ab}+B_{ba}
 =-\frac{1}{2}(g_{il}W^l{}_j+g_{jl}W^l{}_i)A^i{}_aA^j{}_b.
\]
Since \(B_{ab}=B_{ba}\) and \(A\) is invertible, \(B=0\) exactly when \(g_{il}W^l{}_j\) is skew in \(i,j\). This proves part (4).
\end{proof}

If \(\Gamma^i\) has the form~\eqref{eq:quadratic-equation}, then \(V_{ijk}=0\) and, as in the proof of part (3), \(H^i{}_{jk}=\gamma^i{}_{jk}\).
Since \(\gamma^i{}_{jk}\) depends only on \(x\), its derivatives \(\frac{\delta\gamma^i{}_{jk}}{\delta x^l}\) are \(\frac{\partial\gamma^i{}_{jk}}{\partial x^l}\), and by~\eqref{eq:adapted-curvature} the curvature component \(R^i{}_{jkl}\) of the normal Galilean connection is the Riemann curvature tensor of \(g\), with the convention
\[
R\Bigl(\frac{\partial}{\partial x^k},\frac{\partial}{\partial x^l}\Bigr)\frac{\partial}{\partial x^j}=R^i{}_{jkl}\frac{\partial}{\partial x^i},\qquad
R(Y,Z)=\nabla^g_Y\nabla^g_Z-\nabla^g_Z\nabla^g_Y-\nabla^g_{[Y,Z]},
\]
where \(\nabla^g\) is the Levi-Civita connection of \(g\).
With this convention, the sectional curvature of the plane spanned by \(Y\) and \(Z\) is \(K(Y,Z)=\frac{g(R(Y,Z)Z,Y)}{g(Y,Y)g(Z,Z)-g(Y,Z)^2}\).

By Corollary~\ref{cor:local-normal}, there is a normal Galilean connection over each coordinate domain. The connections over two overlapping domains can agree only if their coefficients \(\Gamma^i\) are related by~\eqref{eq:gamma-transition}, which holds exactly when the local vector fields \(\mathscr S_{\mathcal O}\) are restrictions of a single vector field on \(X\).

Every smooth paracompact \(M\) admits a second order differential equation. Choose a vector field \(\mathscr S_\alpha\) of the form \(\mathscr S_{\mathcal O}\) over each coordinate domain \(\mathcal O_\alpha\), and let \((\rho_\alpha)\) be a locally finite partition of unity on \(M\) subordinate to these domains. Then
\[
\mathscr S=\sum_\alpha(\rho_\alpha\circ q)\mathscr S_\alpha
\]
is a smooth vector field on \(X\), with each summand extended by zero outside \(q^{-1}(\mathcal O_\alpha)\). Since \(\mathcal C(\mathscr S_\alpha)=0\), \(\tau(\mathscr S_\alpha)=1\), and \(\sum_\alpha\rho_\alpha=1\), we have \(\mathcal C(\mathscr S)=0\) and \(\tau(\mathscr S)=1\).

\begin{proposition}[The global normal Galilean connection]
\label{prop:global-normal}
Let \((M,g)\) be a smooth paracompact Riemannian manifold and let \(\mathscr S\) be a second order differential equation on \(M\), parameterized by \(t\). On \(P=q^*F(M)\), with the forms \(\tau\) and \(\omega^a\) defined at the beginning of this section, there are regular Galilean connections for which adapted coordinates are laboratory coordinates, the fundamental tensor is \(q^*g\) on \(\mathcal E\), and the geodesics parameterized by \(t\), viewed as curves in \(X\), are exactly the integral curves of \(\mathscr S\). Among these connections there is precisely one normal Galilean connection.
\end{proposition}

\begin{proof}
By Corollary~\ref{cor:local-normal}, over each coordinate domain there is a unique normal Galilean connection with the adapted coordinates as laboratory coordinates, geodesic equations given by \(\mathscr S\), and fundamental tensor \(q^*g\). It remains to check that these local connections agree on overlaps.

Since \(\mathscr S\) is a vector field on \(X\), its coefficients \(\Gamma^i\) satisfy the transformation law~\eqref{eq:gamma-transition} under a change of coordinates \(\bar x=\bar x(x)\) on \(M\). Differentiating that law with respect to \(v\) gives the transformation law for \(N^i{}_j=\frac{1}{2}\frac{\partial\Gamma^i}{\partial v^j}\); together these are the coframe transformations~\eqref{eq:x-coframe-transitions} and~\eqref{eq:v-coframe-transitions}. The orthonormal sections determined by \(g\) differ by an orthogonal matrix \(h\), with
\[
\bar s=sh,\qquad
\bar U^a{}_i\frac{\partial\bar x^i}{\partial x^j}
 =(h^{-1})^a{}_bU^b{}_j.
\]
Consequently, \(\bar s^*\omega=h^{-1}s^*\omega\) and \(\bar s^*\phi=h^{-1}s^*\phi\). Normality is invariant under this change of section. By the uniqueness in Corollary~\ref{cor:local-normal}, the \(\mathfrak{o}_n\)-component \(\Pi\) therefore satisfies the transition laws
\begin{equation}
\begin{split}
\bar s^*\Pi&=h^{-1}(s^*\Pi)h+h^{-1}dh,\\
\bar\Pi^i{}_j\frac{\partial\bar x^j}{\partial x^k}
 &=\frac{\partial\bar x^i}{\partial x^j}\Pi^j{}_k
       -d\left(\frac{\partial\bar x^i}{\partial x^k}\right).
\end{split}
\label{eq:global-connection-overlap}
\end{equation}
These transformations, together with those of \(\tau,\omega,\phi\), are the transition laws of the full Galilean connection. The local connections thus glue on \(P\). Their regularity, fundamental tensor, and geodesic equations agree with the prescribed data over every coordinate domain. Any other normal Galilean connection with those data agrees with this one over each coordinate domain by Corollary~\ref{cor:local-normal}, and hence globally.
\end{proof}

The integral curves of \(\mathscr S\), which need not be complete, are the geodesics of the normal Galilean connection on their domains of definition.

For every smooth paracompact manifold \(M\), \(J^1(\mathbb R,M)\) therefore carries the structure of a regular Galilean manifold: choose a Riemannian metric on \(M\) and, by a partition of unity, a second order differential equation on \(M\).

Conversely, a regular Galilean manifold need not be \(J^1(\mathbb R,M)\) with \(\tau=p^*dt\): in Example~\ref{ex:forced-oscillator} the form \(\tau\) is closed but not exact.

\section{Local and global inverse problems}
\label{sec:inverse-problems}

Fix a Riemannian metric \(g=g(x)\) and a second order differential
equation on \(M\) parameterized by \(t\). By Proposition~\ref{prop:global-normal}, we obtain a unique normal Galilean connection, whose vector field \(\Xi\) is the given second order differential equation. We ask whether \(\Xi\) is the
Euler--Lagrange equation of a regular Lagrangian \(L\) with
\(\frac{\partial^2L}{\partial v^i\partial v^j}=g_{ij}\). The Euler--Lagrange equations of a regular Lagrangian \(L\) are \(\frac{d^2x^i}{dt^2}+\Gamma^i=0\) exactly when
\[
\frac{d}{dt}\frac{\partial L}{\partial v^i}-\frac{\partial L}{\partial x^i}
 =\frac{\partial^2L}{\partial v^i\partial v^j}
  \left(\frac{d^2x^j}{dt^2}+\Gamma^j\right)
\]
along every curve with \(v=\frac{dx}{dt}\); we call \(\frac{\partial^2L}{\partial v^i\partial v^j}\) the multiplier
matrix of \(L\). The general inverse problem, in which the multiplier matrix is unknown, is treated in
\cite{anderson1992inverse,crampin1981differential,crampin1984helmholtz,bucataru2010helmholtz}; here we prescribe it to be \(g_{ij}(x)\).

For a multiplier matrix \(g_{ij}(t,x,v)\), the conditions for \(\frac{d^2x^i}{dt^2}+\Gamma^i=0\) to be the Euler--Lagrange equations of a Lagrangian with multiplier matrix \(g\) near each point are the Helmholtz conditions, given by Douglas \cite[(4.7)--(4.11)]{douglas1941solution}.
In the form of Sarlet \cite[Theorem~1]{sarlet1982helmholtz}, whose matrix \(\Phi^{(0)}\) is \(\mathcal R\), and in the notation of Section~\ref{sec:curvature-bianchi}, they state that \(g\) is symmetric and nondegenerate, that \(\mathring\nabla_\Xi g=0\), that \(g_{ik}\mathcal R^k{}_j\) is symmetric, and that \(\frac{\partial g_{ij}}{\partial v^k}\) is symmetric in all three indices.
Crampin, Prince, and Thompson \cite[Theorem~2]{crampin1984helmholtz} write them as conditions on a two-form on \(J^1(\mathbb R,M)\); see \cite[Theorem~8.5]{sarlet1995derivations} for an intrinsic statement.
In \cite[\S4]{lackey2000metric} they are written as conditions on the torsion of the connection of \cite{lackey2000metric} described in Section~\ref{sec:curvature-bianchi}.
For \(g=g(x)\) the last condition is vacuous, and since \(B\) is half of \(\mathring\nabla_\Xi g\), the remaining conditions are \(B=0\) and the self-adjointness of \(\mathcal R\) (Proposition~\ref{prop:local-inverse}).
Globally, Takens \cite{takens1979global} and Anderson and Duchamp \cite{anderson1980existence} showed that the obstruction to a global Lagrangian for an equation that has Lagrangians near each point is a class in the de Rham cohomology of \(\mathbb R\times M\); see also \cite[\S4.3]{krupkova1997geometry}.
For a prescribed multiplier matrix \(g(x)\), this obstruction vanishes exactly when a closed two-form on \(M\) is exact (Theorem~\ref{thm:global-inverse}).

\begin{corollary}
\label{cor:quadratic-forms}
If \(B=0\) everywhere, the coefficients \(W^i{}_j\) and \(F^i\) in~\eqref{eq:quadratic-equation}
define differential forms on \(M\), depending smoothly on \(t\):
\[
\begin{gathered}
W=\frac{1}{2}W_{ij}\,dx^i\wedge dx^j,\qquad
F=F_i\,dx^i,
\\
W_{ij}=g_{ik}W^k{}_j,\qquad
F_i=g_{ij}F^j.
\end{gathered}
\]
We regard them as forms on \(\mathbb R\times M\) with no \(dt\)
component.
\end{corollary}

\begin{proof}
By Proposition~\ref{prop:riemannian-quadratic}, if \(B=0\) then \(\Gamma^i=\gamma^i{}_{jk}v^jv^k+W^i{}_jv^j+F^i\) with \(g_{il}W^l{}_j\) skew in \(i,j\). Under a change of coordinates on \(M\), \(\Gamma^i\) satisfies the transformation law~\eqref{eq:gamma-transition}, and by the transformation law of the Levi-Civita coefficients, \(\gamma^i{}_{jk}v^jv^k\) has the same inhomogeneous term. Hence \(W^i{}_jv^j+F^i\) transforms as a vector along the projection \(X\to M\). Since the induced change of the coordinates \(v^i\) is linear and independent of \(t\), the coefficients \(W^i{}_j\) and \(F^i\) transform as an endomorphism and a vector on \(M\). With the first index lowered, \(W_{ij}\) is a skew-symmetric two-tensor. Thus \(W\) and \(F\) are globally defined forms on \(M\), depending smoothly on \(t\).
\end{proof}

The forms \(W\) and \(F\) are defined with respect to the identification \(X=J^1(\mathbb R,M)\), whose changes of adapted coordinates do not depend on \(t\). A change of laboratory coordinates that depends on \(t\) changes \(\Gamma^i\), and with it \(W\) and \(F\): such a change removes \(\Gamma\) in Example~\ref{ex:free-fall} and removes \(W\) in Example~\ref{ex:foucault}.

If \(L\) has multiplier matrix \(g\), then
\(L-\frac12g_{ij}v^iv^j\) has zero second derivative in \(v\).
Integrating twice in \(v\) gives
\begin{equation}
L(t,x,v)=\frac{1}{2}g_{ij}(x)v^iv^j+\alpha_i(t,x)v^i+\varphi(t,x).
\label{eq:lagrangian}
\end{equation}
The functions \(\alpha_i\) and \(\varphi\) are the components of the one-form \(\eta=\alpha_i\,dx^i+\varphi\,dt\) on \(\mathbb R\times\mathcal O\), where \(\mathcal O\subset M\) is the domain of the coordinates \(x^i\), and
\[
L=\frac{1}{2}g_{ij}v^iv^j+\eta\Bigl(\frac{\partial}{\partial t}+v^i\frac{\partial}{\partial x^i}\Bigr).
\]
We write \(\alpha=\alpha_i\,dx^i\) for the part of \(\eta\) with no \(dt\) component, so that \(\eta=\alpha+\varphi\,dt\).

\begin{corollary}
\label{cor:euler-lagrange-normal}
If \(\Xi\) is the Euler--Lagrange equation of a Lagrangian with multiplier matrix \(g\), then \(B=0\), and hence \(P^a{}_{bc}=T^a{}_{bc}=0\).
\end{corollary}

\begin{proof}
For~\eqref{eq:lagrangian}, direct calculation gives
\begin{align}
\frac{d}{dt}\frac{\partial L}{\partial v^i}
              -\frac{\partial L}{\partial x^i}
 &=g_{ij}\left(\frac{d^2x^j}{dt^2}+\Gamma^j\right),\nonumber\\
\Gamma^i
 &=\gamma^i{}_{jk}v^jv^k
  +g^{ik}\left(\frac{\partial \alpha_k}{\partial x^j}
                   -\frac{\partial \alpha_j}{\partial x^k}\right)v^j
  +g^{ik}\left(\frac{\partial \alpha_k}{\partial t}
                   -\frac{\partial \varphi}{\partial x^k}\right).
\label{eq:euler-lagrange-gamma}
\end{align}
The coefficient of \(v^j\) in \(\Gamma^i\) is skew with respect to \(g\), so \(B=0\) by Proposition~\ref{prop:riemannian-quadratic}, and then \(P^a{}_{bc}=T^a{}_{bc}=0\) by the same proposition.
\end{proof}

\begin{example}
\label{ex:damped-inverse}
Friction dissipates energy, and motion with friction is the standard example of motion without a Lagrangian of the form kinetic minus potential energy; it is usually described with a dissipation function instead. Nevertheless the damped harmonic oscillator has a Lagrangian, found by Caldirola and Kanai, whose second derivative with respect to \(v\) grows exponentially in \(t\). This example shows that friction rules out a Lagrangian whose multiplier matrix is the Euclidean metric, and that the obstruction disappears when the fundamental tensor is the multiplier matrix of the Caldirola--Kanai Lagrangian. A damped harmonic oscillator, such as a mass on a spring with a viscous damper, obeys
\[
\frac{d^2x}{dt^2}+2\zeta\nu_0\frac{dx}{dt}+\nu_0^2x=0,
\]
where \(\nu_0\) is the natural frequency and \(\zeta\ge0\) the damping ratio \cite[\S25]{landau1976mechanics}. Let \(M=\mathbb R\) with \(g_{11}=1\), so that \(\Gamma=2\zeta\nu_0v+\nu_0^2x\). In the section \(s\) given by \(h=1\), the normal Galilean connection has \(N=\frac{1}{2}\frac{\partial\Gamma}{\partial v}=\zeta\nu_0\) and \(s^*\Pi=0\), and by~\eqref{eq:riemannian-B-coefficients}, \(2B_{11}=-\frac{\partial\Gamma}{\partial v}\), so \(B=-\zeta\nu_0\). Since \(B\ne0\) when \(\zeta\ne0\), Corollary~\ref{cor:euler-lagrange-normal} shows that the equation is not the Euler--Lagrange equation of a Lagrangian with multiplier matrix \(g_{11}=1\). With \(L=\frac{1}{2}v^2-\frac{1}{2}\nu_0^2x^2\) it can instead be written with the Rayleigh dissipation function \(\zeta\nu_0v^2\) as
\[
\frac{d}{dt}\frac{\partial L}{\partial v}
 -\frac{\partial L}{\partial x}
 =-\frac{\partial}{\partial v}\left(\zeta\nu_0v^2\right),
\]
an equation of the type treated in \cite{mestdag2011nonconservative}. It is also the Euler--Lagrange equation of the Lagrangian \(e^{2\zeta\nu_0t}\bigl(\frac{1}{2}v^2-\frac{1}{2}\nu_0^2x^2\bigr)\) of Caldirola and Kanai \cite{caldirola1941forze,kanai1948quantization}, whose multiplier matrix \(e^{2\zeta\nu_0t}\) depends on \(t\). For the fundamental tensor \(g_{11}=e^{2\zeta\nu_0t}\), by~\eqref{eq:B-symmetric} we have \(2B_{11}=\bigl(\Xi(g_{11})-2g_{11}N\bigr)(A^1{}_1)^2=0\). Thus whether \(B\) vanishes depends on the fundamental tensor.
The fundamental tensor \(e^{2\zeta\nu_0t}\) depends on \(t\), which the fundamental theorem allows, while the local and global inverse problems prescribe a multiplier matrix \(g(x)\) on \(M\); the Caldirola--Kanai Lagrangian therefore lies outside them.
\end{example}

By~\eqref{eq:euler-lagrange-gamma}, when \(B=0\) the Lagrangian with one-form \(\eta\) has Euler--Lagrange equations \(\Xi\) if and only if the forms \(W\) and \(F\) satisfy
\begin{equation}
W_{ij}=\frac{\partial \alpha_i}{\partial x^j}
                   -\frac{\partial \alpha_j}{\partial x^i},\qquad
F_i=\frac{\partial \alpha_i}{\partial t}
                   -\frac{\partial \varphi}{\partial x^i}.
\label{eq:coefficient-relations}
\end{equation}
When \(B=0\), define the two-form
\begin{equation}
\mathcal F=-W+dt\wedge F
\quad\text{on }\mathbb R\times M.
\label{eq:two-form-F}
\end{equation}
Writing \(d_M\) for the exterior derivative on \(M\), we have
\[
d\eta=d_M\alpha+dt\wedge\left(\frac{\partial \alpha}{\partial t}-d_M\varphi\right),\qquad
d\mathcal F=-d_M W
             -dt\wedge\left(\frac{\partial W}{\partial t}+d_M F\right).
\]
Hence~\eqref{eq:coefficient-relations} is equivalent to
\begin{equation}
\mathcal F=d\eta.
\label{eq:F-exact}
\end{equation}
If \(\mathcal F=d\eta\), then \(d\mathcal F=0\). Conversely, if \(d\mathcal F=0\) on a contractible open subset of \(\mathbb R\times M\), then by the Poincar\'e lemma \(\mathcal F=d\eta\) there for some one-form \(\eta\).

Assume \(B=0\). In terms of the Levi-Civita connection \(\nabla^g\),
the coordinate form of the Jacobi endomorphism~\eqref{eq:jacobi-endomorphism} is
\begin{align}
 \mathcal R^i{}_j
 &=R^i{}_{kjl}v^kv^l
   +\left(-\frac{1}{2}\nabla^g_kW^i{}_j+\nabla^g_jW^i{}_k\right)v^k - \frac{1}{2}\frac{\partial W^i{}_j}{\partial t}
       -\frac{1}{4}W^i{}_kW^k{}_j+\nabla^g_jF^i.
\label{eq:quadratic-jacobi-endomorphism}
\end{align}
Here \(R^i{}_{kjl}\) are the components of the Riemann curvature tensor of \(g\).
After the first index is lowered, the curvature term and \(W^i{}_kW^k{}_j\) are symmetric. Since \(\nabla^g\) preserves \(g\) and \(W\) is skew, the skew part of \(\mathcal R\) is
\begin{align}
 \mathcal R_{ij}-\mathcal R_{ji}
 &= -\left[\left(\frac{\partial W_{ij}}{\partial x^k}
      +\frac{\partial W_{ki}}{\partial x^j}
      +\frac{\partial W_{jk}}{\partial x^i}\right)v^k
      +\frac{\partial W_{ij}}{\partial t}
      -\frac{\partial F_i}{\partial x^j}
      +\frac{\partial F_j}{\partial x^i}\right].
\label{eq:skew-jacobi-tensor}
\end{align}
\begin{proposition}[The local inverse problem]
\label{prop:local-inverse}
Suppose \(g_{ij}(x)\) and \(\Gamma^i(t,x,v)\) are given.
The differential equations are locally Euler--Lagrange equations
with multiplier matrix \(g\) if and only if their normal Galilean connection
satisfies
\[
B=0,\qquad \mathcal R_{ab}=\mathcal R_{ba}
\quad(\text{equivalently } \mathcal R_{ij}=\mathcal R_{ji}).
\]
When \(B=0\), \(\mathcal R\) is self-adjoint with respect to \(g\) if and only if
\begin{equation}
d_M W=0,\qquad
\frac{\partial W}{\partial t}+d_M F=0,
\label{eq:local-closure}
\end{equation}
that is, if and only if \(d\mathcal F=0\).
\end{proposition}

\begin{proof}
By Corollary~\ref{cor:euler-lagrange-normal}, \(B=0\) is necessary; assume it. By~\eqref{eq:skew-jacobi-tensor}, \(\mathcal R_{ij}-\mathcal R_{ji}\) is affine in \(v\). Its coefficient of \(v^k\) vanishes for all \(i,j,k\) exactly when \(d_MW=0\), and its term independent of \(v\) vanishes exactly when \(\frac{\partial W}{\partial t}+d_MF=0\). Thus \(\mathcal R\) is self-adjoint if and only if~\eqref{eq:local-closure} holds. Since \(d\mathcal F=-d_MW-dt\wedge\bigl(\frac{\partial W}{\partial t}+d_MF\bigr)\), this is equivalent to \(d\mathcal F=0\).

If \(\Xi\) is the Euler--Lagrange equation of a Lagrangian with multiplier matrix \(g\) on an open set, then, since~\eqref{eq:coefficient-relations} is equivalent to~\eqref{eq:F-exact}, \(\mathcal F=d\eta\) there for the one-form \(\eta\) of that Lagrangian, so \(d\mathcal F=0\) and \(\mathcal R\) is self-adjoint. Conversely, suppose \(d\mathcal F=0\), and let \(\mathcal O\subset M\) be a contractible coordinate neighborhood and \(I\) an interval. By the Poincar\'e lemma there is a one-form \(\eta\) on \(I\times\mathcal O\) with \(\mathcal F=d\eta\), and the Lagrangian \(\frac12g_{ij}v^iv^j+\eta\bigl(\frac{\partial}{\partial t}+v^i\frac{\partial}{\partial x^i}\bigr)\) has Euler--Lagrange equations \(\Xi\) over \(I\times\mathcal O\).
\end{proof}

\begin{example}
\label{ex:foucault}
In a rotating frame of reference a moving body is deflected by the Coriolis force, which is perpendicular to its velocity and proportional to it. Since this force depends on the velocity, it is not the gradient of a potential, and the motion has no Lagrangian of the form kinetic minus potential energy. Nevertheless the motion has a Lagrangian with multiplier matrix \(g\), every such Lagrangian containing a term linear in \(v\). The Coriolis force is also removed by passing to a suitably rotating frame, which is a change of laboratory coordinates depending on \(t\). This example exhibits both phenomena. The Foucault pendulum at latitude \(\vartheta\) is a pendulum of natural frequency \(\nu_0\) observed from the Earth, which rotates with angular speed \(\nu_E\). For small oscillations, neglecting the horizontal component of the Earth's angular velocity and absorbing the centrifugal force into the effective gravity, the horizontal displacement \(x\) of the bob, in coordinates fixed to the Earth with \(x^1\) pointing east and \(x^2\) north, satisfies
\[
\frac{d^2x}{dt^2}+2\nu_z\mathcal J\frac{dx}{dt}+\nu_0^2x=0,\qquad
\nu_z=\nu_E\sin\vartheta,\qquad
\mathcal J=\begin{pmatrix}0&-1\\1&0\end{pmatrix},
\]
where \(-2\nu_z\mathcal J\frac{dx}{dt}\) is the Coriolis force per unit mass \cite[\S39]{landau1976mechanics}. Let \(M=\mathbb R^2\) with \(g_{ij}=\delta_{ij}\). Here \(\Gamma^i=W^i{}_jv^j+F^i\) with \(W^i{}_j=2\nu_z\mathcal J^i{}_j\), which is skew, and \(F^i=\nu_0^2x^i\); thus \(B=0\) by Proposition~\ref{prop:riemannian-quadratic}, \(W=-2\nu_z\,dx^1\wedge dx^2\) is constant, and \(F=d\bigl(\frac{1}{2}\nu_0^2\lvert x\rvert^2\bigr)\). Hence \(d_MW=0\) and \(\frac{\partial W}{\partial t}+d_MF=0\), and by Proposition~\ref{prop:local-inverse} the equation is an Euler--Lagrange equation with multiplier matrix \(\delta_{ij}\). Indeed, the Lagrangian
\[
L=\frac{1}{2}\lvert v\rvert^2+\nu_z(x^1v^2-x^2v^1)-\frac{1}{2}\nu_0^2\lvert x\rvert^2
\]
has \(\alpha=\nu_z(x^1\,dx^2-x^2\,dx^1)\) and \(\varphi=-\frac{1}{2}\nu_0^2\lvert x\rvert^2\), which satisfy~\eqref{eq:coefficient-relations}. The one-form \(\alpha\) of any Lagrangian with multiplier matrix \(\delta_{ij}\) for this equation satisfies \(d_M\alpha=-W=2\nu_z\,dx^1\wedge dx^2\ne0\), since \(\mathcal F=d\eta\): the Coriolis force is what requires \(\alpha\).

Introduce rotating coordinates \(y\) by \(x=e^{-\nu_zt\mathcal J}y\), where \(e^{-\nu_zt\mathcal J}\) is the rotation by the angle \(-\nu_zt\). Since \(\frac{d}{dt}e^{-\nu_zt\mathcal J}=-\nu_z\mathcal Je^{-\nu_zt\mathcal J}\) and \(\mathcal J^2=-\mathbf 1\),
\[
\frac{d^2x}{dt^2}+2\nu_z\mathcal J\frac{dx}{dt}+\nu_0^2x
 =e^{-\nu_zt\mathcal J}\Bigl(\frac{d^2y}{dt^2}+(\nu_0^2+\nu_z^2)\,y\Bigr).
\]
In the coordinates \(y\), the coefficient \(W\) of the equation is zero.
\end{example}

The global inverse problem asks for a one-form \(\eta\) with \(\mathcal F=d\eta\) on all of \(\mathbb R\times M\). This condition reduces to the exactness of a two-form on \(M\).
\begin{theorem}[The global inverse problem]
\label{thm:global-inverse}
Let \((M,g)\) be a smooth paracompact Riemannian manifold, and let
\(\Xi\) be the vector field of the normal Galilean connection determined
by \(g\) and a second order differential equation on \(M\),
parameterized by \(t\). The following
are equivalent:
\begin{enumerate}
\item There is a smooth global Lagrangian with multiplier matrix \(g\)
whose Euler--Lagrange equations are those of \(\Xi\).
\item \(B=0\), and \(\mathcal F=d\eta\) for a one-form \(\eta\) on \(\mathbb R\times M\), where \(\mathcal F=-W+dt\wedge F\).
\item \(B=0\), \(\mathcal R_{ab}=\mathcal R_{ba}\), and the two-form \(W(0)\) on \(M\), the restriction of \(W\) to \(\{0\}\times M\), is exact.
\end{enumerate}
\end{theorem}

\begin{proof}
\((1)\Rightarrow(2)\). Let \(L\) be a global Lagrangian with multiplier matrix \(g\) whose Euler--Lagrange equations are those of \(\Xi\). In adapted coordinates \(L\) has the form~\eqref{eq:lagrangian}. A change of adapted coordinates acts linearly on the coordinates \(v^i\), so it preserves the terms of \(L\) of each degree in \(v\); hence the one-form \(\eta\) is defined on all of \(\mathbb R\times M\). By Corollary~\ref{cor:euler-lagrange-normal}, \(B=0\), and then \(\mathcal F=d\eta\) because \(\eta\) satisfies~\eqref{eq:coefficient-relations}.

\((2)\Rightarrow(1)\). If \(B=0\) and \(\mathcal F=d\eta\) for a one-form \(\eta\) on \(\mathbb R\times M\), then the Lagrangian \(L=\frac12g_{ij}v^iv^j+\eta\bigl(\frac{\partial}{\partial t}+v^i\frac{\partial}{\partial x^i}\bigr)\) is defined on all of \(X\). Since \(\mathcal F=d\eta\) is equivalent to~\eqref{eq:coefficient-relations}, its Euler--Lagrange equations are those of \(\Xi\).

\((2)\Rightarrow(3)\). If \(\mathcal F=d\eta\), then \(d\mathcal F=0\), so \(\mathcal R\) is self-adjoint by Proposition~\ref{prop:local-inverse}. The restrictions of \(\mathcal F\) and \(\eta\) to \(\{0\}\times M\) are \(-W(0)\) and \(\alpha(0)\), so \(W(0)=-d_M\alpha(0)\).

\((3)\Rightarrow(2)\). By Proposition~\ref{prop:local-inverse}, \(d_MW=0\) and \(\frac{\partial W}{\partial t}+d_MF=0\). Choose a one-form \(\alpha_0\) on \(M\) with \(d_M\alpha_0=-W(0)\), and set
\begin{equation}
\alpha(t)=\alpha_0+\int_0^t F(s)\,ds,\qquad \varphi=0.
\label{eq:global-temporal-gauge}
\end{equation}
Then
\[
d_M\alpha(t)=-W(0)+\int_0^t d_M F(s)\,ds
 =-W(0)-\int_0^t\frac{\partial W}{\partial s}(s)\,ds=-W(t),\qquad
\frac{\partial \alpha}{\partial t}=F,
\]
so \(\eta=\alpha\) satisfies \(d\eta=-W+dt\wedge F=\mathcal F\).
\end{proof}

Since \(d_MW=0\), the two-form \(W(0)\) is closed. Thus, once \(B=0\) and \(\mathcal R\) is self-adjoint, the existence of a global Lagrangian with multiplier matrix \(g\) depends only on the cohomology class of \(W(0)\) in \(H^2_{\mathrm{dR}}(M)\). In particular, if \(H^2_{\mathrm{dR}}(M)=0\), every equation that is locally an Euler--Lagrange equation with multiplier matrix \(g\) is the Euler--Lagrange equation of a global Lagrangian with multiplier matrix \(g\). On \(M=S^2\) this fails (Example~\ref{ex:monopole}).

The one-form \(\eta\) in condition (2) is determined up to a closed one-form. If \(\eta=\alpha+\varphi\,dt\) and \(\eta'=\alpha'+\varphi'\,dt\) both satisfy \(\mathcal F=d\eta\), then \(\beta=\eta'-\eta\) is a closed one-form on \(\mathbb R\times M\), and the corresponding Lagrangians differ by
\[
\beta\Bigl(\frac{\partial}{\partial t}+v^i\frac{\partial}{\partial x^i}\Bigr)
 =(\alpha'_i-\alpha_i)v^i+\varphi'-\varphi.
\]
For this function, \(\frac{d}{dt}\frac{\partial}{\partial v^i}-\frac{\partial}{\partial x^i}\) gives \(d\beta\bigl(\frac{\partial}{\partial t}+v^j\frac{\partial}{\partial x^j},\frac{\partial}{\partial x^i}\bigr)=0\), so its Euler--Lagrange equations hold identically: the difference is a null Lagrangian, and every null Lagrangian has this form \cite[Corollary~4.3.8]{krupkova1997geometry}. If \(\beta=df\) for a function \(f\) on \(\mathbb R\times M\), the difference is the total time derivative \(\frac{df(t,x(t))}{dt}\). Since \(H^1_{\mathrm{dR}}(\mathbb R\times M)\cong H^1_{\mathrm{dR}}(M)\), when \(H^1_{\mathrm{dR}}(M)=0\) any two global Lagrangians with multiplier matrix \(g\) for the same equation differ by a total time derivative \cite[Theorem~4.3]{anderson1980existence}. On \(M=S^1\) this fails (Example~\ref{ex:flux-ring}).

\begin{example}
\label{ex:flux-ring}
The simplest configuration space containing a closed curve that bounds no surface is the circle. On the circle a force can do nonzero work around the closed curve, as the electric field induced by a changing magnetic flux does. Such a force has no potential, so the motion it produces has no global Lagrangian of the form kinetic minus potential energy. Nevertheless the motion can have a global Lagrangian with multiplier matrix \(g\), whose term linear in \(v\) depends on \(t\). The circle also carries closed one-forms that are not exact, and these give null Lagrangians that are not total time derivatives, so two global Lagrangians for the same equation need not differ by a total time derivative. This example exhibits both phenomena. A bead with charge-to-mass ratio \(\varepsilon\) slides on a circular ring of radius \(r\) threaded by a long solenoid that carries the magnetic flux \(\Psi(t)\). The magnetic field vanishes on the ring, but by Faraday's law a changing flux induces an electric field along it. Let \(M=S^1\) with the metric \(g=r^2\kappa\otimes\kappa\), where \(\kappa\) is the closed one-form with \(\int_{S^1}\kappa=2\pi\), which is \(d\theta\) in every angular coordinate \(\theta\). Per unit mass, the Lagrangian of the bead is
\[
L=\frac{1}{2}r^2\Bigl(\frac{d\theta}{dt}\Bigr)^2+\frac{\varepsilon\Psi(t)}{2\pi}\frac{d\theta}{dt}.
\]
Its one-form is \(\eta=\frac{\varepsilon\Psi(t)}{2\pi}\kappa\), so \(W=-d_M\alpha=0\) and \(F=\frac{\partial\alpha}{\partial t}=\frac{\varepsilon}{2\pi}\frac{d\Psi}{dt}\kappa\). The equation of motion is \(r^2\frac{d^2\theta}{dt^2}=-\frac{\varepsilon}{2\pi}\frac{d\Psi}{dt}\). Since \(\int_{S^1}F=\varepsilon\frac{d\Psi}{dt}\) is, up to sign, the induced electromotive force per unit mass, \(F\) is closed but not exact whenever \(\frac{d\Psi}{dt}\ne0\). Nevertheless \(W(0)=0\) is exact, so Theorem~\ref{thm:global-inverse} gives a global Lagrangian; taking \(\alpha_0=\frac{\varepsilon\Psi(0)}{2\pi}\kappa\) in~\eqref{eq:global-temporal-gauge} recovers \(\eta\). Since \(L\) does not depend on \(\theta\), the canonical angular momentum \(r^2\frac{d\theta}{dt}+\frac{\varepsilon\Psi}{2\pi}\) per unit mass is conserved: changing the flux from \(\Psi_1\) to \(\Psi_2\) changes \(r^2\frac{d\theta}{dt}\) by \(-\frac{\varepsilon}{2\pi}(\Psi_2-\Psi_1)\), however the change is made.

If \(\Psi\) is constant, then \(F=0\), and the equation \(r^2\frac{d^2\theta}{dt^2}=0\) is also the Euler--Lagrange equation of \(\frac{1}{2}r^2\bigl(\frac{d\theta}{dt}\bigr)^2\). The two Lagrangians differ by the null Lagrangian \(\frac{\varepsilon\Psi}{2\pi}\frac{d\theta}{dt}\), which is not a total time derivative because \(\kappa\) is not exact. This is the classical counterpart of the configuration of Aharonov and Bohm \cite{aharonov1959significance}: the constant flux has no effect on the classical motion, but it changes the phase of the quantum mechanical wave function.
\end{example}

\begin{example}
\label{ex:monopole}
A magnetic monopole has nonzero total magnetic flux through any sphere centered on it. Near every point of such a sphere, the magnetic force on a charged particle confined to the sphere is described by a Lagrangian with a term linear in the velocity, but because the total flux is nonzero there is no such Lagrangian on the whole sphere. This example shows that the motion has a Lagrangian with multiplier matrix \(g\) near every point but none on the whole sphere, and gives the local Lagrangians explicitly. A charged particle constrained to a sphere centered on a magnetic monopole \cite{dirac1931quantised,haldane1983fractional} moves in a magnetic field normal to the sphere and of constant strength. With the radius of the sphere normalized to \(1\), let \(M=S^2\) with its standard metric and area form \(\mathrm{vol}_{S^2}\). The equation of motion is \(\frac{d^2x^i}{dt^2}+\gamma^i{}_{jk}v^jv^k+W^i{}_jv^j=0\), with \(\gamma^i{}_{jk}\) the Levi-Civita coefficients and \(W=\nu_c\,\mathrm{vol}_{S^2}\), where \(W^i{}_j\) is obtained by raising the first index of \(W\) and the constant \(\nu_c>0\), the cyclotron frequency, is the product of the charge-to-mass ratio and the field strength. Since \(W\) has top degree, \(d_MW=0\), and \(\frac{\partial W}{\partial t}+d_MF=0\); thus~\eqref{eq:local-closure} holds, and by Proposition~\ref{prop:local-inverse} the equation is locally an Euler--Lagrange equation with multiplier matrix \(g\). However, \(W(0)=\nu_c\,\mathrm{vol}_{S^2}\) is not exact: by Stokes' theorem an exact two-form on the closed surface \(S^2\) has integral zero, whereas \(\int_{S^2}W(0)=4\pi\nu_c\). By Theorem~\ref{thm:global-inverse}, there is no global Lagrangian with multiplier matrix \(g\).

In colatitude \(\theta\) and longitude \(\lambda\), let \(\alpha_N\) and \(\alpha_S\) be \(-\nu_c\) times the unit-strength monopole potentials of \cite[Section III, equation (11)]{wu1975nonintegrable}:
\[
\alpha_N=\nu_c(\cos\theta-1)\,d\lambda,\qquad
\alpha_S=\nu_c(\cos\theta+1)\,d\lambda.
\]
Although \(d\lambda\) is singular at both poles, in the coordinates \((x^1,x^2,x^3)\) of \(\mathbb R^3\) we have \(x^1\,dx^2-x^2\,dx^1=\sin^2\theta\,d\lambda\), so
\[
\alpha_N=-\nu_c\frac{x^1\,dx^2-x^2\,dx^1}{1+x^3},\qquad
\alpha_S=\nu_c\frac{x^1\,dx^2-x^2\,dx^1}{1-x^3}.
\]
Thus \(\alpha_N\) is smooth on \(\mathcal O_N=S^2\setminus\{\text{south pole}\}\) and \(\alpha_S\) on \(\mathcal O_S=S^2\setminus\{\text{north pole}\}\). Each satisfies \(d_M\alpha=-W\) on its domain, so with \(\varphi=0\) the Lagrangians \(\frac{1}{2}g_{ij}v^iv^j+(\alpha_N)_iv^i\) and \(\frac{1}{2}g_{ij}v^iv^j+(\alpha_S)_iv^i\) have the given equation on \(\mathcal O_N\) and \(\mathcal O_S\), respectively. On the overlap, \(\alpha_S-\alpha_N=2\nu_c\,d\lambda\) is closed but not exact, so these Lagrangians differ by the null Lagrangian \(2\nu_c\frac{d\lambda}{dt}\), which is a total time derivative only locally.
\end{example}

\section{Variations and the Jacobi equation}
\label{sec:jacobi}

Let \((M,g)\) be a Riemannian manifold, \(X=J^1(\mathbb R,M)\) with adapted coordinates \((t,x^i,v^i)\), and \(P=q^*F(M)\). A second order differential equation on \(M\) parameterized by \(t\) determines a normal Galilean connection on \(P\) whose vector field \(\Xi\) on \(X\) is that equation (Proposition~\ref{prop:global-normal}). We ask how the integral curves of \(\Xi\) near a given integral curve deviate from it. To first order, this deviation satisfies a linear second order equation along the curve, the Jacobi equation. When \(\Xi\) is the geodesic equation of \(g\), the endomorphism in the Jacobi equation is given by the Riemann curvature tensor of \(g\). For a general second order differential equation it is the Jacobi endomorphism \(\mathcal R\) \cite{crampin1996linear}, and the Jacobi equation is studied in the geometric theory of second order differential equations \cite{antonelli2001invariants,bucataru2011geometric,sarlet1995derivations}. The eigenvalues of \(\mathcal R\) are used as a criterion for the stability of nearby integral curves, called Jacobi stability \cite{bohmer2012jacobi,punzi2009geometry}.

The Jacobi equation obtained below from the structure equations of the normal Galilean connection is the classical one: since \(\nabla^{(t)}-B\) is the dynamical covariant derivative \(\mathring\nabla_\Xi\) (Section~\ref{sec:curvature-bianchi}), it is the equation \(\mathring\nabla_\Xi^2\xi+\mathcal R\xi=0\) of \cite[Theorem~4]{crampin1996linear}.
The normal Galilean connection splits \(\mathring\nabla_\Xi\) as \(\nabla^{(t)}-B\), where \(\nabla^{(t)}\) preserves \(g\) and \(B\) is half the dynamical covariant derivative of \(g\).
When \(B=0\) and \(\mathcal R\) is self-adjoint, which are the conditions of the local inverse problem, the metric \(g\) therefore bounds the Jacobi fields, provided \(\nabla^{(t)}\mathcal R\) has integrable norm along the integral curve (Theorem~\ref{thm:jacobi-bound}).

Let \(\sigma:I\to X\) be an integral curve of \(\Xi\), where \(I\) is an interval. A \emph{variation} of \(\sigma\) is a smooth map
\[
\rho:I\times(-u_0,u_0)\longrightarrow X,\qquad
\rho(t,u)=\bigl(t,x^i(t,u),v^i(t,u)\bigr)
\]
in adapted coordinates, such that \(\rho(t,0)=\sigma(t)\) and each curve \(t\mapsto\rho(t,u)\) is an integral curve of \(\Xi\); thus all of these curves share the time parameter \(t\), and \(v^i=\frac{\partial x^i}{\partial t}\). Its \emph{variation field} along \(\sigma\) has components \(y^i=\frac{\partial x^i}{\partial u}\) at \(u=0\).

Let \(s\) be a local section of \(P\); shrinking \(I\) and \(u_0\) if necessary, the image of \(\rho\) lies in its domain. Since \(\rho\) preserves \(t\), \(\rho^*(s^*\tau)=dt\). Since each curve \(t\mapsto\rho(t,u)\) is a geodesic, the proof of Proposition~\ref{prop:geodesic-equations} shows that \(\rho^*(s^*\omega^a)\) and \(\rho^*(s^*\phi^a)\) vanish on \(\frac{\partial}{\partial t}\). Hence
\[
\rho^*(s^*\omega^a)=\xi^a\,du,\qquad
\rho^*(s^*\phi^a)=\psi^a\,du
\]
for functions \(\xi^a,\psi^a\) on \(I\times(-u_0,u_0)\). In adapted coordinates, \(s^*\omega^a=(A^{-1})^a{}_i(dx^i-v^i\,dt)\), so \(\xi^a=(A^{-1})^a{}_i\frac{\partial x^i}{\partial u}\) is the variation field in the orthonormal frame of \(\mathcal E\) dual to \(s^*\omega^a\). Since \(\rho_*\bigl(\frac{\partial}{\partial t}\bigr)=\Xi\), we have \(\rho^*(s^*\Pi^a{}_b)=(s^*\Pi^a{}_b)(\Xi)\,dt\) modulo \(du\).

We pull back the definitions~\eqref{eq:spatial-torsion} and~\eqref{eq:velocity-torsion} of \(T^a\) and \(\Phi^a\) by \(s\circ\rho\). In the torsion decompositions of Theorem~\ref{thm:metric-general}(1), every term other than those with \(\tau\wedge\omega^b\) and \(\tau\wedge\phi^b\) pulls back to a multiple of \(du\wedge du=0\), and \(D=0\) by normality. Thus
\[
\rho^*(s^*T^a)=B^a{}_b\xi^b\,dt\wedge du,\qquad
\rho^*(s^*\Phi^a)=\bigl(-\mathcal R^a{}_b\xi^b+B^a{}_b\psi^b\bigr)\,dt\wedge du,
\]
while
\[
\rho^*\bigl(s^*(d\omega^a+\Pi^a{}_b\wedge\omega^b+\phi^a\wedge\tau)\bigr)
 =\Bigl(\frac{\partial\xi^a}{\partial t}+(s^*\Pi^a{}_b)(\Xi)\xi^b-\psi^a\Bigr)dt\wedge du,
\]
\[
\rho^*\bigl(s^*(d\phi^a+\Pi^a{}_b\wedge\phi^b)\bigr)
 =\Bigl(\frac{\partial\psi^a}{\partial t}+(s^*\Pi^a{}_b)(\Xi)\psi^b\Bigr)dt\wedge du.
\]
Evaluating the formula \((\nabla\xi)^a=d\xi^a+\Pi^a{}_b\xi^b\) on \(\Xi\) along each curve \(t\mapsto\rho(t,u)\) gives the temporal covariant derivative
\[
(\nabla^{(t)}\xi)^a=\frac{\partial\xi^a}{\partial t}+(s^*\Pi^a{}_b)(\Xi)\xi^b.
\]
Comparing coefficients of \(dt\wedge du\) therefore gives
\begin{align}
\psi&=\nabla^{(t)}\xi-B\xi,\label{eq:variation-psi}\\
\nabla^{(t)}\psi&=-\mathcal R\xi+B\psi.\label{eq:variation-psi-derivative}
\end{align}

\begin{proposition}
\label{prop:variation-equation}
For a variation \(\rho\) of an integral curve \(\sigma\) of \(\Xi\), the components \(\xi^a\) satisfy
\begin{equation}
(\nabla^{(t)})^2\xi-2B\nabla^{(t)}\xi
 =\bigl(-\mathcal R+\nabla^{(t)}B-B^2\bigr)\xi.
\label{eq:variation-second-order}
\end{equation}
\end{proposition}

\begin{proof}
Applying \(\nabla^{(t)}\) to~\eqref{eq:variation-psi} and using~\eqref{eq:variation-psi-derivative} gives
\[
(\nabla^{(t)})^2\xi-(\nabla^{(t)}B)\xi-B\nabla^{(t)}\xi
 =-\mathcal R\xi+B(\nabla^{(t)}\xi-B\xi).
\]
Rearranging terms then gives~\eqref{eq:variation-second-order}.
\end{proof}

The equations \(\psi=\nabla^{(t)}\xi-B\xi\) and \(\nabla^{(t)}\psi=-\mathcal R\xi+B\psi\) can be written
using the derivative \(\nabla^{(t)}-B\):
\[
\psi=(\nabla^{(t)}-B)\xi,\qquad
(\nabla^{(t)}-B)\psi=-\mathcal R\xi.
\]
Consequently, the Jacobi equation is
\begin{equation}
(\nabla^{(t)}-B)^2\xi=-\mathcal R\xi,
\label{eq:canonical-jacobi}
\end{equation}
where \(\mathcal R\) is the Jacobi endomorphism. A \emph{Jacobi field} along \(\sigma\) is a solution \(\xi\) of~\eqref{eq:canonical-jacobi} along \(\sigma\); by Proposition~\ref{prop:variation-equation}, the variation field of every variation of \(\sigma\) is a Jacobi field. The derivative
\(\nabla^{(t)}-B\) agrees with \(\nabla^{(t)}\) precisely when \(B=0\); Example~\ref{ex:damped} shows the effect of \(B\ne0\).

For the geodesic equation of \(g\), with \(\Gamma^i=\gamma^i{}_{jk}v^jv^k\), we have \(N^i{}_j=\gamma^i{}_{jk}v^k\), and \(B=0\) by Proposition~\ref{prop:riemannian-quadratic}; thus \(\nabla^{(t)}\) is the covariant derivative of the Levi-Civita connection along \(\sigma\). By~\eqref{eq:quadratic-jacobi-endomorphism} with \(W=0\) and \(F=0\), \(\mathcal R^i{}_j=R^i{}_{kjl}v^kv^l\), where \(R^i{}_{kjl}\) is the Riemann curvature tensor of \(g\). Thus the Jacobi equation is the classical one,
\begin{equation}
\bigl((\nabla^{(t)})^2\xi\bigr)^i
   +R^i{}_{kjl}\,\frac{dx^k}{dt}\frac{dx^l}{dt}\xi^j=0.
\label{eq:riemannian-jacobi}
\end{equation}
For an equation other than the geodesic equation, \(\mathcal R\) need not vanish when \(g\) is flat (Example~\ref{ex:tides}).

In coordinates, the Jacobi equation is the linearization of the differential equation along \(\sigma\). Differentiating \(\frac{\partial^2x^i}{\partial t^2}+\Gamma^i\bigl(t,x,\frac{\partial x}{\partial t}\bigr)=0\) with respect to \(u\) at \(u=0\) shows that the variation field satisfies
\begin{equation}
\frac{d^2y^i}{dt^2}+\frac{\partial\Gamma^i}{\partial v^j}\frac{dy^j}{dt}
            +\frac{\partial\Gamma^i}{\partial x^j}y^j=0.
\label{eq:coordinate-linearization}
\end{equation}
In the coordinate frame along \(\sigma\), by~\eqref{eq:nabla-temporal},
\begin{align*}
\nabla^{(t)}-B&=\frac{d}{dt}+N,\qquad
N^i{}_j=\frac{1}{2}\frac{\partial\Gamma^i}{\partial v^j},\\
\mathcal R^i{}_j&=\frac{\partial\Gamma^i}{\partial x^j}
          -\frac{dN^i{}_j}{dt}-N^i{}_kN^k{}_j.
\end{align*}
Since \(\rho^*(s^*\omega^a)=\xi^a\,du\) and \(\rho^*(s^*\phi^a)=\psi^a\,du\), at \(u=0\) we have
\(y=A\xi\) and \(A\psi=\frac{dy}{dt}+Ny\). Thus
\((\nabla^{(t)}-B)\) acts on \(y\) as \(\frac{d}{dt}+N\), and
\[
\left(\frac{d}{dt}+N\right)^2y+\mathcal R y
 =\frac{d^2y}{dt^2}+2N\frac{dy}{dt}
   +\frac{\partial\Gamma}{\partial x}y.
\]
Since \(2N^i{}_j=\frac{\partial\Gamma^i}{\partial v^j}\), this is exactly
\eqref{eq:coordinate-linearization}.

\begin{example}
\label{ex:tides}
In general relativity, the relative acceleration of two nearby freely falling bodies is given by the equation of geodesic deviation, whose endomorphism is defined by the Riemann curvature tensor of spacetime \cite[Chapter 11]{misner1973gravitation}. This example shows that for Newtonian gravity the Jacobi equation of the normal Galilean connection is the Newtonian equation of geodesic deviation, with the tidal tensor in place of the Riemann curvature tensor. In Newtonian gravity, a body moving in the gravitational potential \(\mathcal V\) of a mass distribution of density \(\varrho\) obeys
\[
\frac{d^2x^i}{dt^2}+\delta^{ij}\frac{\partial\mathcal V}{\partial x^j}=0,\qquad
\delta^{ij}\frac{\partial^2\mathcal V}{\partial x^i\partial x^j}=4\pi G\varrho,
\]
on \(M=\mathbb R^3\) with \(g_{ij}=\delta_{ij}\), where \(G\) is the gravitational constant and the gravitational force per unit mass is \(-\delta^{ij}\frac{\partial\mathcal V}{\partial x^j}\). Since \(\Gamma^i\) does not depend on \(v\), we have \(W=0\), so \(B=0\) by Proposition~\ref{prop:riemannian-quadratic}. In the section \(s\) given by \(h=1\), \(N=0\) and \(s^*\Pi=0\), so \(\nabla^{(t)}=\frac{d}{dt}\) and
\[
\mathcal R_{ij}=\frac{\partial^2\mathcal V}{\partial x^i\partial x^j},
\]
the tidal tensor. The Jacobi equation \(\frac{d^2\xi}{dt^2}+\mathcal R\xi=0\) is the Newtonian equation of geodesic deviation, which governs the relative acceleration of two nearby freely falling bodies \cite[Chapter 12]{misner1973gravitation}, and \(\operatorname{tr}\mathcal R=4\pi G\varrho\). Outside a point mass \(m_0\), \(\mathcal V=-\frac{Gm_0}{r}\) with \(r=\lvert x\rvert\), and
\[
\mathcal R=\frac{Gm_0}{r^3}\bigl(\mathbf 1-3\hat x\hat x^{T}\bigr),\qquad \hat x=\frac{x}{r},
\]
with the eigenvalue \(-\frac{2Gm_0}{r^3}\) in the radial direction and \(\frac{Gm_0}{r^3}\) on its orthogonal complement: nearby freely falling bodies separate radially and approach transversally, which is the origin of the tides. The curvature \(R^i{}_{kjl}\) of \(g\) vanishes, while \(\mathcal R\ne0\); the gravitational field appears in the Jacobi endomorphism of the normal Galilean connection, as it appears in the curvature of the spacetime connection of Newton--Cartan geometry \cite{kunzle1972galilei,misner1973gravitation}.
\end{example}

\begin{example}
\label{ex:damped}
Friction makes nearby motions of a damped oscillator approach one another at a rate set by the damping, whether the motion is underdamped, critically damped, or overdamped. This example shows that in the Jacobi equation of the normal Galilean connection these two effects separate: the approach is produced by the torsion coefficient \(B\) alone, and the Jacobi endomorphism \(\mathcal R\) distinguishes the three regimes. This continues Example~\ref{ex:damped-inverse}, where the damped harmonic oscillator is defined by
\[
\frac{d^2x}{dt^2}+2\zeta\nu_0\frac{dx}{dt}+\nu_0^2x=0,
\]
with natural frequency \(\nu_0\) and damping ratio \(\zeta\ge0\). This is a second order differential equation on \(M=\mathbb R\), for which we take \(g_{11}=1\). In the section \(s\) given by \(h=1\), its normal Galilean connection has
\[
N=\zeta\nu_0,\qquad s^*\Pi=0,\qquad B=-\zeta\nu_0,\qquad \mathcal R=\nu_0^2(1-\zeta^2),
\]
where \(\mathcal R\) follows from~\eqref{eq:jacobi-endomorphism}. Thus \(\nabla^{(t)}=\frac{d}{dt}\) and \(\nabla^{(t)}-B=\frac{d}{dt}+\zeta\nu_0\), so the solutions of the Jacobi equation \((\nabla^{(t)}-B)^2\xi=-\mathcal R\xi\) are \(\xi=e^{-\zeta\nu_0t}f\) with \(\frac{d^2f}{dt^2}+\mathcal Rf=0\). The sign of \(\mathcal R\) separates the three regimes of the damped oscillator: underdamped motion with frequency \(\sqrt{\mathcal R}=\nu_0\sqrt{1-\zeta^2}\) when \(\zeta<1\), critical damping when \(\zeta=1\), and overdamped motion when \(\zeta>1\). The decay factor \(e^{-\zeta\nu_0t}\) is produced by \(B\) alone.
\end{example}

Jacobi stability concerns whether integral curves of \(\Xi\) near \(\sigma\) remain near it, that is, whether Jacobi fields along \(\sigma\) remain bounded. Consider first the geodesic equation of \(g\). There \(B=0\), and \(\mathcal R\) is self-adjoint with respect to \(g\), with
\[
g(\mathcal Rz,z)=K(z,v)\bigl(g(v,v)\,g(z,z)-g(z,v)^2\bigr),
\]
where \(K(z,v)\) is the sectional curvature of the plane spanned by \(z\) and \(v\). Since \(\mathcal Rv=0\), Jacobi fields tangent to \(\sigma\) have the form \((a+bt)v\) and arise from reparameterizing \(\sigma\); stability concerns Jacobi fields orthogonal to \(v\) \cite[Chapter 5]{docarmo1992riemannian}. Their behavior is governed by the sign of the sectional curvature \cite[Chapters 5, 7, 9 and 10]{docarmo1992riemannian}, \cite[Chapter 1]{cheeger1975comparison}. If \(K\le0\), then
\[
\frac{d^2}{dt^2}\,\frac{1}{2}g(\xi,\xi)=g(\nabla^{(t)}\xi,\nabla^{(t)}\xi)-g(\mathcal R\xi,\xi)\ge0,
\]
so a nonzero Jacobi field vanishing at one time grows at least linearly and has no further zeros; nearby geodesics diverge, which underlies the Cartan--Hadamard theorem. If \(K\ge\kappa>0\), the Rauch comparison theorem shows that, along a geodesic with \(g(v,v)=1\), every Jacobi field vanishing at one time vanishes again within time \(\frac{\pi}{\sqrt\kappa}\), which underlies the Bonnet--Myers theorem. Thus under positive curvature orthogonal Jacobi fields oscillate, but positive curvature by itself does not bound them (see the surface following Example~\ref{ex:mathieu}).

For a general second order differential equation and vector fields \(\xi_1,\xi_2\) along \(\sigma\), since \(\nabla^{(t)}\) preserves \(g\) and \(B_{ab}=B_{ba}\) by normality,
\[
\frac{d}{dt}\,g(\xi_1,\xi_2)
 =g\bigl((\nabla^{(t)}-B)\xi_1,\xi_2\bigr)+g\bigl(\xi_1,(\nabla^{(t)}-B)\xi_2\bigr)+2g(B\xi_1,\xi_2),
\]
so the derivative \(\nabla^{(t)}-B\) of the Jacobi equation preserves \(g\) exactly when \(B=0\). In that case the Jacobi equation is \((\nabla^{(t)})^2\xi=-\mathcal R\xi\). By the solution of the local inverse problem, the conditions \(B=0\) and \(\mathcal R\) self-adjoint with respect to \(g\) hold exactly when \(\Xi\) is locally the Euler--Lagrange equation of a Lagrangian with multiplier matrix \(g\).
The Jacobi equation is then equivalent to the Euler--Lagrange equation of the second variation of such a Lagrangian, which is the Jacobi equation of the calculus of variations \cite[Chapter~5, \S\S27 and~29]{gelfand1963calculus}.

\begin{theorem}
\label{thm:jacobi-bound}
Let \(\sigma:I\to X\) be an integral curve of \(\Xi\), with \(t_0\in I\). Along \(\sigma\), suppose that \(B=0\), that \(\mathcal R\) is self-adjoint with respect to \(g\), and that for some \(\mu>0\)
\[
g(\mathcal Rz,z)\ge\mu\,g(z,z),\qquad t\in I.
\]
Let \(\|\nabla^{(t)}\mathcal R\|\) be the operator norm with respect to \(g\). For every Jacobi field \(\xi\) along \(\sigma\), the function
\[
E=g(\nabla^{(t)}\xi,\nabla^{(t)}\xi)+g(\mathcal R\xi,\xi)
\]
satisfies
\[
\mu\,g(\xi,\xi)\le E,\qquad g(\nabla^{(t)}\xi,\nabla^{(t)}\xi)\le E,\qquad
E(t)\le E(t_0)\exp\Bigl(\frac{1}{\mu}\Bigl|\int_{t_0}^t\|\nabla^{(t)}\mathcal R\|\,ds\Bigr|\Bigr).
\]
In particular, if \(\int_I\|\nabla^{(t)}\mathcal R\|\,dt<\infty\), then every Jacobi field along \(\sigma\) is bounded on \(I\).
\end{theorem}

\begin{proof}
Since \(B=0\), the Jacobi equation is \((\nabla^{(t)})^2\xi=-\mathcal R\xi\). Since \(\nabla^{(t)}\) preserves \(g\), and \(\mathcal R\) and \(\nabla^{(t)}\mathcal R\) are self-adjoint,
\[
\frac{dE}{dt}
 =2g\bigl((\nabla^{(t)})^2\xi,\nabla^{(t)}\xi\bigr)+g\bigl((\nabla^{(t)}\mathcal R)\xi,\xi\bigr)+2g(\mathcal R\xi,\nabla^{(t)}\xi)
 =g\bigl((\nabla^{(t)}\mathcal R)\xi,\xi\bigr).
\]
The lower bound on \(\mathcal R\) gives \(\mu\,g(\xi,\xi)\le g(\mathcal R\xi,\xi)\le E\) and \(g(\nabla^{(t)}\xi,\nabla^{(t)}\xi)\le E\). Hence
\[
\Bigl|\frac{dE}{dt}\Bigr|\le\|\nabla^{(t)}\mathcal R\|\,g(\xi,\xi)\le\frac{1}{\mu}\|\nabla^{(t)}\mathcal R\|\,E,
\]
and Gr\"onwall's inequality gives the bound on \(E\).
\end{proof}

The proof uses the self-adjointness of \(\mathcal R\) to cancel the terms of \(\frac{dE}{dt}\) involving \(\nabla^{(t)}\xi\), and the integrability of \(\|\nabla^{(t)}\mathcal R\|\) to bound \(E\) on all of \(I\). Neither hypothesis can be omitted (Examples~\ref{ex:optical-trap} and~\ref{ex:mathieu}).

When \(B=0\) and \(\mathcal R\) is self-adjoint, in a frame along \(\sigma\) that is parallel for \(\nabla^{(t)}\), the Jacobi equation is the Euler--Lagrange equation of \(\ell=\frac12g\bigl(\frac{d\xi}{dt},\frac{d\xi}{dt}\bigr)-\frac12g(\mathcal R\xi,\xi)\), and \(E\) is twice its energy \(g\bigl(\frac{d\xi}{dt},\frac{d\xi}{dt}\bigr)-\ell\); if \(\nabla^{(t)}\mathcal R=0\), then \(\ell\) does not depend on \(t\), and \(E\) is constant. For a scalar equation \(\frac{d^2f}{dt^2}+\mathcal R(t)f=0\) with \(\mathcal R\ge\mu>0\), we have \(\|\nabla^{(t)}\mathcal R\|=\bigl|\frac{d\mathcal R}{dt}\bigr|\), and the theorem is the statement that all solutions are bounded when \(\mathcal R\) has bounded variation \cite[Chapter 6, Section 4]{bellman1953stability}.

\begin{corollary}
\label{cor:riemannian-jacobi-bound}
Let \(\Xi\) be the geodesic equation of \(g\), let \(c:I\to M\) be a nonconstant geodesic of \(g\), and let \(\sigma=\bigl(t,c,\frac{dc}{dt}\bigr)\) be the corresponding integral curve of \(\Xi\), with \(v=\frac{dc}{dt}\). Suppose that \(K(z,v)\ge\kappa>0\) for every \(z\) not parallel to \(v\), and that \(\int_I\|\nabla^{(t)}\mathcal R\|\,dt<\infty\). Then every Jacobi field along \(\sigma\) orthogonal to \(v\) is bounded on \(I\).
\end{corollary}

\begin{proof}
For the geodesic equation, \(B=0\) by Proposition~\ref{prop:riemannian-quadratic}, and \(\mathcal R\) is self-adjoint by the symmetries of the Riemann curvature tensor. Since \(v\) is parallel, \(\nabla^{(t)}\) preserves the orthogonal complement of \(v\), and \(g(v,v)\) is a positive constant. Since \(\mathcal R\) is self-adjoint and \(\mathcal Rv=0\), it also preserves this complement, on which \(g(\mathcal Rz,z)\ge\kappa\,g(v,v)\,g(z,z)\). The proof of Theorem~\ref{thm:jacobi-bound}, applied to Jacobi fields with values in this complement and with \(\mu=\kappa\,g(v,v)\), gives the result.
\end{proof}

Along a geodesic, \(\nabla^{(t)}\mathcal R\) is determined by the covariant derivative of the curvature tensor of \(g\) in the direction \(v\). If \(g\) is locally symmetric, then \(\nabla^{(t)}\mathcal R=0\); thus on a locally symmetric space with \(K\ge\kappa>0\), such as a round sphere, Jacobi fields orthogonal to a geodesic are bounded.

\begin{example}
\label{ex:optical-trap}
Near a stable equilibrium, a restoring force makes nearby motions oscillate about one another, and they remain close. A circulatory force, which acts perpendicular to the displacement and is not the gradient of a potential, leaves the restoring force unchanged but can nevertheless destabilize the equilibrium. In the Jacobi equation the restoring force gives the self-adjoint part of \(\mathcal R\), and the circulatory force its skew-adjoint part. This example shows, for a particle in an optical trap, that the skew-adjoint part makes the Jacobi fields grow exponentially although the self-adjoint part is a positive constant. A dielectric particle held in vacuum by a circularly polarized Gaussian laser beam experiences, near the beam axis, a restoring optical force together with an azimuthal force due to the spin of the light, in an approximately constant ratio to the restoring force; at low pressure this azimuthal force ejects the particle from the trap \cite{svak2018transverse}. In the plane transverse to the beam, with \(M=\mathbb R^2\) and \(g_{ij}=\delta_{ij}\), and neglecting gas damping, the linearized motion is
\[
\frac{d^2x}{dt^2}+\nu_0^2(\mathbf 1+\varepsilon\mathcal J)x=0,\qquad
\mathcal J=\begin{pmatrix}0&-1\\1&0\end{pmatrix},
\]
where \(\nu_0\) is the trap frequency and \(\varepsilon>0\) is the ratio of the azimuthal to the restoring force. Since \(\Gamma^i\) does not depend on \(v\), we have \(W=0\), so \(B=0\) by Proposition~\ref{prop:riemannian-quadratic}. In the section \(s\) given by \(h=1\), \(s^*\Pi=0\), so \(\nabla^{(t)}=\frac{d}{dt}\) and \(\mathcal R=\nu_0^2(\mathbf 1+\varepsilon\mathcal J)\). Thus \(\nabla^{(t)}\mathcal R=0\) and, since \(\mathcal J\) is skew, \(g(\mathcal Rz,z)=\nu_0^2g(z,z)\); but \(\mathcal R\) is not self-adjoint, and, as a point of interest, the equation is not locally the Euler--Lagrange equation of a Lagrangian with multiplier matrix \(\delta_{ij}\), by Proposition~\ref{prop:local-inverse}. Putting
\[
\vartheta=\nu_0\sqrt{\frac{\sqrt{1+\varepsilon^2}-1}{2}},\qquad
\Theta=\nu_0\sqrt{\frac{\sqrt{1+\varepsilon^2}+1}{2}},
\]
the Jacobi field
\[
\xi(t)=e^{\vartheta t}\begin{pmatrix}\cos\Theta t\\-\sin\Theta t\end{pmatrix}
\]
satisfies \(g(\xi,\xi)=e^{2\vartheta t}\): the particle spirals out of the trap at the rate \(\vartheta\). With gas damping \(2\zeta\nu_0\frac{dx}{dt}\), for which \(B\ne0\), the trap is stable if and only if \(2\zeta>\varepsilon\), which gives the pressure below which the particle escapes.
\end{example}

\begin{example}
\label{ex:mathieu}
A pendulum whose restoring force is modulated periodically in time can be driven to large amplitude, although the restoring force never falls below a positive fraction of its mean; this is parametric resonance. In the Jacobi equation the modulated restoring force gives a Jacobi endomorphism \(\mathcal R\) that is bounded below by a positive constant but varies periodically. This example shows that such an \(\mathcal R\) has unbounded Jacobi fields, so a positive lower bound on \(\mathcal R\) does not bound the Jacobi fields unless the variation of \(\mathcal R\) is controlled. A pendulum of natural frequency \(\nu_0\) whose point of suspension is driven vertically at frequency \(2\nu_0\), with amplitude small compared with its length, obeys, for small angles \(x\), the Mathieu equation
\[
\frac{d^2x}{dt^2}+\nu_0^2(1+\varepsilon\cos2\nu_0t)\,x=0,\qquad 0<\varepsilon<1,
\]
where \(\varepsilon\) is proportional to the driving amplitude \cite[\S27]{landau1976mechanics}. Let \(M=\mathbb R\) with \(g_{11}=1\). Since \(\Gamma^i\) does not depend on \(v\), we have \(W=0\), so \(B=0\) by Proposition~\ref{prop:riemannian-quadratic}, and in the section \(s\) given by \(h=1\), \(\nabla^{(t)}=\frac{d}{dt}\). Here \(\mathcal R=\nu_0^2(1+\varepsilon\cos2\nu_0t)\ge\nu_0^2(1-\varepsilon)>0\) is self-adjoint, but \(\nabla^{(t)}\mathcal R=\frac{d\mathcal R}{dt}\) is periodic and not identically zero, so \(\int_{\mathbb R}\bigl|\frac{d\mathcal R}{dt}\bigr|\,dt=\infty\). For sufficiently small \(\varepsilon\), the driving frequency \(2\nu_0\) lies in the first instability region of the Mathieu equation, so the equation has unbounded solutions \cite[\S27]{landau1976mechanics}; these are unbounded Jacobi fields. This is parametric resonance: to first order in \(\varepsilon\), the unbounded solutions grow like \(e^{\varepsilon\nu_0t/4}\) \cite[\S27]{landau1976mechanics}. For \(\nu_0=1\) and \(\varepsilon=\frac25\), the monodromy matrix over one period \(\pi\) has trace \(-2.099\), and the growth rate is \(0.0995\), against \(\frac14\varepsilon\nu_0=0.1\).
\end{example}

The same unbounded Jacobi fields occur along a geodesic of a surface whose sectional curvature is bounded below by a positive constant. Let \(K(t)=1+\varepsilon\cos2t\) with \(0<\varepsilon<1\) as in Example~\ref{ex:mathieu}, with \(\nu_0=1\). On the strip \(\mathbb R\times(-1,1)\) with coordinates \((t,w)\), the metric \(\bigl(1-\frac{1}{2}K(t)w^2\bigr)^2dt^2+dw^2\) is defined, since \(K<2\). It has the geodesic \(w=0\), parameterized by arc length \(t\), along which the sectional curvature is \(K(t)\ge1-\varepsilon\). Its Jacobi fields \(f(t)\frac{\partial}{\partial w}\) orthogonal to this geodesic satisfy \(\frac{d^2f}{dt^2}+K(t)f=0\), and some of them are unbounded. Thus a lower bound on the sectional curvature by a positive constant does not by itself bound the Jacobi fields orthogonal to a geodesic.

\section*{Acknowledgments}

We would like to thank James Grant and Paul Busch for useful discussions concerning the direction of this paper.

\section*{AI Disclosure}

We used large language models (Claude Opus and ChatGPT) in preparing this paper: to draft and revise text and examples, and to check computations symbolically.
The author verified all mathematical content and takes full responsibility for this work.

\bibliographystyle{plain}
\bibliography{references}

\begin{thebibliography}{10}

\bibitem{aharonov1959significance}
Y.~Aharonov and D.~Bohm.
\newblock Significance of electromagnetic potentials in the quantum theory.
\newblock {\em Physical Review}, 115(3):485--491, 1959.

\bibitem{anderson1992inverse}
Ian Anderson and Gerard Thompson.
\newblock The inverse problem of the calculus of variations for ordinary
  differential equations.
\newblock {\em Memoirs of the American Mathematical Society}, 98(473), 1992.

\bibitem{anderson1980existence}
Ian~M. Anderson and T.~Duchamp.
\newblock On the existence of global variational principles.
\newblock {\em American Journal of Mathematics}, 102(5):781--868, 1980.

\bibitem{andringa2011newtonian}
Roel Andringa, Eric Bergshoeff, Sudhakar Panda, and Mees de~Roo.
\newblock Newtonian gravity and the {Bargmann} algebra.
\newblock {\em Classical and Quantum Gravity}, 28(10):105011, 2011.

\bibitem{antonelli2003kcc}
P.~L. Antonelli and I.~Bucataru.
\newblock {KCC}-theory of a system of second order differential equations.
\newblock In Peter~L. Antonelli, editor, {\em Handbook of {Finsler} Geometry},
  volume~1, pages 83--174. Kluwer Academic Publishers, Dordrecht, 2003.

\bibitem{antonelli2001invariants}
Peter~L. Antonelli and Ioan Buc{\u{a}}taru.
\newblock New results about the geometric invariants in {KCC}-theory.
\newblock {\em Analele {\c{S}}tiin{\c{t}}ifice ale Universit{\u{a}}{\c{t}}ii
  ``Al. I. Cuza'' din Ia{\c{s}}i. Matematic{\u{a}}}, 47(2):405--420, 2001.

\bibitem{bao2012introduction}
David Bao, Shiing-Shen Chern, and Zhongmin Shen.
\newblock {\em An introduction to Riemann-Finsler geometry}, volume 200 of {\em
  Graduate Texts in Mathematics}.
\newblock Springer, New York, 2000.

\bibitem{bellman1953stability}
Richard Bellman.
\newblock {\em Stability Theory of Differential Equations}.
\newblock McGraw-Hill, New York, 1953.

\bibitem{bohmer2012jacobi}
Christian~G. B{\"o}hmer, Tiberiu Harko, and Sorin~V. Sabau.
\newblock Jacobi stability analysis of dynamical systems---applications in
  gravitation and cosmology.
\newblock {\em Advances in Theoretical and Mathematical Physics},
  16(4):1145--1196, 2012.

\bibitem{bucataru2007metric}
Ioan Bucataru.
\newblock Metric nonlinear connections.
\newblock {\em Differential Geometry and its Applications}, 25(3):335--343,
  2007.

\bibitem{bucataru2010helmholtz}
Ioan Bucataru and Oana Constantinescu.
\newblock {Helmholtz} conditions and symmetries for the time dependent case of
  the inverse problem of the calculus of variations.
\newblock {\em Journal of Geometry and Physics}, 60(11):1710--1725, 2010.

\bibitem{bucataru2011geometric}
Ioan Bucataru, Oana Constantinescu, and Matias~F. Dahl.
\newblock A geometric setting for systems of ordinary differential equations.
\newblock {\em International Journal of Geometric Methods in Modern Physics},
  8(6):1291--1327, 2011.

\bibitem{caldirola1941forze}
P.~Caldirola.
\newblock Forze non conservative nella meccanica quantistica.
\newblock {\em Il Nuovo Cimento}, 18:393--400, 1941.

\bibitem{cheeger1975comparison}
Jeff Cheeger and David~G. Ebin.
\newblock {\em Comparison Theorems in Riemannian Geometry}, volume~9 of {\em
  North-Holland Mathematical Library}.
\newblock North-Holland, Amsterdam, 1975.

\bibitem{chern1939geometrie}
Shiing-Shen Chern.
\newblock Sur la g{\'e}om{\'e}trie d'un syst{\`e}me d'{\'e}quations
  diff{\'e}rentielles du second ordre.
\newblock {\em Bull. Sci. Math.}, 63:206--212, 1939.

\bibitem{crampin1981differential}
Michael Crampin.
\newblock On the differential geometry of the {E}uler-{L}agrange equations, and
  the inverse problem of lagrangian dynamics.
\newblock {\em Journal of Physics A: Mathematical and General},
  14(10):2567--2575, 1981.

\bibitem{crampin1996linear}
Michael Crampin, Eduardo Mart{\'i}nez, and Willy Sarlet.
\newblock Linear connections for systems of second-order ordinary differential
  equations.
\newblock {\em Annales de l'Institut Henri Poincar{\'e}. Physique
  th{\'e}orique}, 65(2):223--249, 1996.

\bibitem{crampin1984helmholtz}
Michael Crampin, Geoff~E. Prince, and Gerard Thompson.
\newblock A geometrical version of the {Helmholtz} conditions in time-dependent
  {Lagrangian} dynamics.
\newblock {\em Journal of Physics A: Mathematical and General},
  17(7):1437--1447, 1984.

\bibitem{dirac1931quantised}
P.~A.~M. Dirac.
\newblock Quantised singularities in the electromagnetic field.
\newblock {\em Proceedings of the Royal Society of London. Series A},
  133(821):60--72, 1931.

\bibitem{docarmo1992riemannian}
Manfredo~Perdig{\~a}o do~Carmo.
\newblock {\em Riemannian Geometry}.
\newblock Mathematics: Theory \& Applications. Birkh{\"a}user, Boston, 1992.

\bibitem{doubrov1999equivalence}
Boris Doubrov, Boris Komrakov, and Tohru Morimoto.
\newblock Equivalence of holonomic differential equations.
\newblock {\em Lobachevskii J. Math}, 3:39--71, 1999.

\bibitem{douglas1941solution}
Jesse Douglas.
\newblock Solution of the inverse problem of the calculus of variations.
\newblock {\em Transactions of the American Mathematical Society},
  50(1):71--128, 1941.

\bibitem{duval1985bargmann}
C.~Duval, G.~Burdet, H.~P. K{\"u}nzle, and M.~Perrin.
\newblock {Bargmann} structures and {Newton-Cartan} theory.
\newblock {\em Physical Review D}, 31(8):1841--1853, 1985.

\bibitem{fels1995equivalence}
M.~E. Fels.
\newblock The equivalence problem for systems of second-order ordinary
  differential equations.
\newblock {\em Proceedings of the London Mathematical Society}, 71(1):221--240,
  1995.

\bibitem{gelfand1963calculus}
I.~M. Gelfand and S.~V. Fomin.
\newblock {\em Calculus of Variations}.
\newblock Prentice-Hall, Englewood Cliffs, NJ, 1963.
\newblock Translated from the Russian by R. A. Silverman.

\bibitem{grant2012galilean}
James D.~E. Grant and Bradley~C. Lackey.
\newblock On {G}alilean connections and the first jet bundle.
\newblock {\em Central European Journal of Mathematics}, 10(5):1889--1895,
  2012.

\bibitem{haldane1983fractional}
F.~D.~M. Haldane.
\newblock Fractional quantization of the {H}all effect: a hierarchy of
  incompressible quantum fluid states.
\newblock {\em Physical Review Letters}, 51(7):605--608, 1983.

\bibitem{kanai1948quantization}
E.~Kanai.
\newblock On the quantization of the dissipative systems.
\newblock {\em Progress of Theoretical Physics}, 3(4):440--442, 1948.

\bibitem{krupkova1997geometry}
Olga Krupkov{\'a}.
\newblock {\em The geometry of ordinary variational equations}, volume 1678 of
  {\em Lecture Notes in Mathematics}.
\newblock Springer, Berlin, 1997.

\bibitem{kunzle1972galilei}
Hans~Peter K{\"u}nzle.
\newblock {Galilei} and {Lorentz} structures on space-time: comparison of the
  corresponding geometry and physics.
\newblock {\em Annales de l'Institut Henri Poincar{\'e}. Physique
  th{\'e}orique}, 17(4):337--362, 1972.

\bibitem{lackey2000metric}
Brad Lackey.
\newblock Metric equivalence of path spaces.
\newblock {\em Nonlinear Studies}, 7(2):241--250, 2002.

\bibitem{landau1975fields}
L.~D. Landau and E.~M. Lifshitz.
\newblock {\em The Classical Theory of Fields}, volume~2 of {\em Course of
  Theoretical Physics}.
\newblock Pergamon Press, Oxford, 4 edition, 1975.

\bibitem{landau1976mechanics}
L.~D. Landau and E.~M. Lifshitz.
\newblock {\em Mechanics}, volume~1 of {\em Course of Theoretical Physics}.
\newblock Pergamon Press, Oxford, 3 edition, 1976.

\bibitem{massa1994jet}
Enrico Massa and Enrico Pagani.
\newblock Jet bundle geometry, dynamical connections, and the inverse problem
  of {Lagrangian} mechanics.
\newblock {\em Annales de l'Institut Henri Poincar{\'e}. Physique
  th{\'e}orique}, 61(1):17--62, 1994.

\bibitem{mestdag2011nonconservative}
Tom Mestdag, Willy Sarlet, and Michael Crampin.
\newblock The inverse problem for {Lagrangian} systems with certain
  non-conservative forces.
\newblock {\em Differential Geometry and its Applications}, 29(1):55--72, 2011.

\bibitem{miron2003geometry}
Radu Miron, Mihai Anastasiei, and I.~Bucataru.
\newblock The geometry of {L}agrange spaces.
\newblock In Peter~L. Antonelli, editor, {\em Handbook of {Finsler} Geometry},
  volume~2, pages 969--1122. Kluwer Academic Publishers, Dordrecht, 2003.

\bibitem{miron2001geometry}
Radu Miron, Drago{\c{s}} Hrimiuc, Hideo Shimada, and Sorin~V. Sab{\u{a}}u.
\newblock {\em The geometry of {H}amilton and {L}agrange spaces}, volume 118 of
  {\em Fundamental Theories of Physics}.
\newblock Kluwer Academic Publishers, Dordrecht, 2001.

\bibitem{misner1973gravitation}
Charles~W. Misner, Kip~S. Thorne, and John~Archibald Wheeler.
\newblock {\em Gravitation}.
\newblock W. H. Freeman, San Francisco, 1973.

\bibitem{punzi2009geometry}
Raffaele Punzi and Mattias N.~R. Wohlfarth.
\newblock Geometry and stability of dynamical systems.
\newblock {\em Physical Review E}, 79(4):046606, 2009.

\bibitem{sarlet1982helmholtz}
W.~Sarlet.
\newblock The {Helmholtz} conditions revisited. {A} new approach to the inverse
  problem of {Lagrangian} dynamics.
\newblock {\em Journal of Physics A: Mathematical and General},
  15(5):1503--1517, 1982.

\bibitem{sarlet1995derivations}
W.~Sarlet, A.~Vandecasteele, F.~Cantrijn, and E.~Mart{\'i}nez.
\newblock Derivations of forms along a map: the framework for time-dependent
  second-order equations.
\newblock {\em Differential Geometry and its Applications}, 5(2):171--203,
  1995.

\bibitem{sharpe1997differential}
Richard~W. Sharpe.
\newblock {\em Differential geometry: {Cartan}'s generalization of {Klein}'s
  {Erlangen} program}, volume 166 of {\em Graduate Texts in Mathematics}.
\newblock Springer, New York, 1997.

\bibitem{souriau1997structure}
Jean-Marie Souriau.
\newblock {\em Structure of dynamical systems: a symplectic view of physics},
  volume 149 of {\em Progress in Mathematics}.
\newblock Birkh{\"a}user, Boston, 1997.

\bibitem{svak2018transverse}
Vojt{\v{e}}ch Svak, Oto Brzobohat{\'y}, Martin {\v{S}}iler, Petr J{\'a}kl, Jan
  Ka{\v{n}}ka, Pavel Zem{\'a}nek, and Stephen~H. Simpson.
\newblock Transverse spin forces and non-equilibrium particle dynamics in a
  circularly polarized vacuum optical trap.
\newblock {\em Nature Communications}, 9, 2018.

\bibitem{takens1979global}
Floris Takens.
\newblock A global version of the inverse problem of the calculus of
  variations.
\newblock {\em Journal of Differential Geometry}, 14(4):543--562, 1979.

\bibitem{wu1975nonintegrable}
Tai~Tsun Wu and Chen~Ning Yang.
\newblock Concept of nonintegrable phase factors and global formulation of
  gauge fields.
\newblock {\em Physical Review D}, 12(12):3845--3857, 1975.

\end{thebibliography}

\end{document}